\magnification=\magstep 1
\input amstex
\documentstyle{amsppt}

\def\a{\alpha}
\def\b{\beta}

\def\O{\Omega}

\def\ti{\widetilde}

\def\e{\epsilon}

\def\d{\delta}
\def\A{{\Cal A}}
\def\pa{\partial}
\def\np{{\hskip 0.8pt \hbox to 0pt{/\hss}\hskip -0.8pt\partial}}
\def\G{{\Cal G}}

\def\B{{\Cal B}}
\def\F{{\Cal F}}

\def\DD{{\Cal D}}
\def\M{{\Cal M}}
\def\N{{\Cal N}}
\def\NN{{\bold N}}
\def\MM{\bold M}
\def\p{\phi}
\def\P{\Phi}
\def\R{{\Cal R}}

\def\n{\nabla}

\def\lmd{\lambda}

\def\Ga{\Gamma}

\def\X{{Cal X}}
\def\S{\Cal S}
\def\SS{\bold  S}
\def\ss{\bold s}

\def\NN{\bold N}

\def\GG{\bold G}
\def\gg{\bold g}
\def\W{{\Cal W}}
\def\SW{\hbox{\bf SW}}
\def\sw{\hbox{\bf sw}}
\def\cs{\hbox{\bf cs}}

\def\AtiyahBott{[2]}

\def\AustinBraamI{[4]}
\def\AustinBraamII{[5]}
\def\BGV{[6]}
\def\Borel{[7]}
\def\Cartan{[8]}
\def\Chen{[9]}

\def\DonaldsonII{[12]}

\def\FloerII{[14]}

\def\Lim{[18]}
\def\MarcolliWang{[19]}

\def\WangYeI{[23]}
\def\WangYe{24}

\def\Ye{[26]}
\def\YeGromov{[27]}
\def\X{{\Cal X}}

\topmatter
\title   Equivariant and Bott-type Seiberg-Witten
Floer homology: Part II
\endtitle
\author Rugang Ye
\endauthor
\address Department of Mathematics, University of
California, Santa Barbara, CA 93106, USA
\endaddress
\email yer\@math.ucsb.edu
\endemail
\address Ruhr-Universit\"at Bochum, Fakult\"at f\"ur Mathematik,
44780 Bochum, Germany
\endaddress
\email ye\@dgeo.ruhr-uni-bochum.de
\endemail

%\keywords Seiberg-Witten equation, Seiberg-Witten Floer homology, Bott
%type 
%\endkeywords
%\subjclass 58
%\endsubjclass
\rightheadtext{Equivariant  Seiberg-Witten
Floer homology } 
\abstract
We construct equivariant and Bott-type 
Seiberg-Witten Floer homology 
and cohomology for 3-manifolds, 
in particular rational homology spheres, and prove their
diffeomorphism invariance.
\endabstract
\endtopmatter
\document

\head Table of Contents
\endhead

\noindent 1. Introduction

\noindent 2. Singular version of  equivariant 
theory

\noindent 3. Invariance

\noindent 4. de Rham version

\noindent 5. Cartan version

\noindent 6. Exponential  convergence

\noindent 7. Structure near infinity

\bigskip

\head{1. Introduction}
\endhead

This paper is the sequel of [$\WangYe$], and Part II of our series of 
papers on Bott-type and equivariant Seiberg-Witten Floer homology.
As explained in the introduction to Part I (see [$\WangYe$]),
our goal is to construct 
diffeomorphism invariants for 3-manifolds, in particular 
rational homology spheres,
which are based on the Seiberg-Witten theory 
and the Floer homology theory. The major problem which we have to 
resolve in the process of constructions is that of
 reducible Seiberg-Witten points
and reducible transition trajectories.
In [$\WangYe$], we constructed 
the Bott-type Seiberg-Witten  Floer homology 
and cohomology and provided  part of the 
proof for their invariance, which is based on the
key technique 
``spinor perturbation". (The remaining part will be 
presented in $\Ye$.)  
Besides, a number of  basic ingredients of 
analyitical or geometrical nature for our theory  
were presented. 

In the present paper, we construct equivariant 
Seiberg-Witten Floer homologies and 
cohomologies and prove their invariance. 
A fundamental ingredient in our constructions (both 
in Part I and the present Part II) is a theory 
on the Morse-Floer-Bott
 flow complex along with its projection 
to critical submanifolds.
Part of this theory  was presented in
Part I. 
In this paper, we complete the theory. 
Most of these results were
 obtained when $\WangYeI$, 
the earlier  version of [$\WangYe$], was written in 1996. 
To arrange more balanced  size of our papers, 
we place them in the present Part II rather than 
in Part I.

We continue with the set-up, notations and terminologies 
in [$\WangYe$].  In particular,
we work in the framework of the Seiberg-Witten theory 
on a given rational homology sphere $Y$.  The following is 
a brief account of the above topics.

Our basic equations are the (perturbed)
3-dimensional Seiberg-Witten
equation    $\sw_{\lambda, H}(a, \p)=0$ on 
$Y$, i.e. 
$$
\cases
*_YF_a+\langle e_i\cdot\p,\p\rangle e^i& = \nabla H(a),\cr
\np_a\p + \lambda \p & =0 \cr
\endcases
\tag 1.1
$$
with perturbations $H$ (holonomy) and 
$\lambda$ (a real number), and the (perturbed) 
Seiberg-Witten trajectory equation (flow equation) 
$\SW_{\lambda, 
H}(A, \P) =0$ on $X=Y \times \Bbb R$, i.e.
$$
\cases
\frac{\partial a}{\partial t}-*_YF_a-d_Yf-\langle e_i\cdot\p,\p\rangle e^i
&=\n H(a), \cr
\frac{\partial \p}{\partial t}+ \np_a\p+\lmd\p+f\p&=0\cr
\endcases
\tag 1.2
$$
 with $A= a +f dt, \P= \p$,
 cf. [$\WangYe$]. 
 Our basic geometrical objects, on which we 
 build our theories, are the 
 moduli space $\R^0$ of based gauge
 classes of Seiberg-Witten points
 (solutions of (1.1))
 and the  moduli spaces of Seiberg-Witten 
 trajectories or flow lines (solutions of (1.2)).
   The fundamental functional involved 
 is the Seiberg-Witten type Chern-Simons functional $\cs$ on the 
 configuration space $\A(Y) \times \Ga(Y)$ and 
 its based gauge quotient $\B^0=
 \A(Y) \times \Ga(Y)) \slash
 \G^0$.
 
\bigskip
\noindent {\bf Equivariant Constructions: Singular Version}
\smallskip

The basic idea is to couple 
the configuration space $\A(Y) \times \Ga(Y)$
  with a space $\SS$ which has 
a free $S^1$ action.  If we first divide out by
the based gauge group $\G^0$, then what we do is to 
couple the based gauge quotient $\B^0 $ with 
$\SS$.  We have the diagonal action of $S^1$ on the 
product $\B^0 \times \SS$. 
  The  
Chern-Simons functional 
is  extended to  the product $\B^0 
\times \SS$ in the trivial way, namely $\cs(\a, \ss) = \cs(\a)$. It is 
invariant under the $S^1$ action,
 and hence  
descends to the quotient 
$\B^0 \times_{S^1} \SS$. 
Its critical submanifolds are given by the 
$S^1$ 
quotient  $\R^0 \times_{S^1} \SS$ of the 
{\it stablized moduli space} $\R^0 \times \SS$.  We use 
{\it generalized cubical singular chains and cochains} on them
to build our basic chains and cochains.  
The boundary operator will then be constructed by 
utilizing the classical boundary operator 
along with {\it stablized Morse-Floer-Bott flow complex}, i.e. 
the stablized  
moduli spaces of (based gauge classes of ) 
Seiberg-Witten 
trajectories,
or rather the $S^1$
quotient of this complex, which we 
call the {\it quotient stablized Morse-Floer-Bott
flow complex}.   Here, the stablization again 
means multiplying by
$\SS$. 
For each choice of 
$\SS$, we obtain in this fashion an
equivariant Seiberg-Witten Floer homology and 
cohomology.

The most obvious choices for $\SS$ are the 
odd dimensional spheres $S^{2n+1}$, $n=0, 1, ...,$ and the 
infinite dimensional sphere $S^{\infty}$.
Indeed, these  choices are natural from the viewpoint of 
the classical 
equivariant homology and cohomology construction due to Borel $\Borel$, and the viewpoint 
 of equivariant Morse theory. (But we were not led 
 to these choices this way. Instead, we arrived at them 
 in an intuitive way.)
Recall that for a topological space $\X$ with 
an action by
a  group $G$,
its $G$-equivariant homology and cohomology as introduced 
by Borel
$\Borel$ are defined to be the  homology and cohomology of
the homotopy quotient space $EG \times_{G} \X$, where $EG$
 denotes
a contractible space with a free $G$ action (thus $BG= EG \slash
G$ is the classifying space of $G$), and $G$ acts on 
$EG \times \X$
by the diagonal action.   Now, if $\B^0$ were a
compact manifold  and the Chern-Simons  functional were 
an ordinary  equivariant Morse function, then our
construction with  $\SS= S^{\infty}= ES^1$ would reproduce
the classical equivariant  homology and cohomology of $\B^0$,
namely the homology and cohomology of the homotopy
quotient $S^{\infty} \times_{S^1} \B^0$.   Of course, the goal
of our constructions is precisely  to produce new invariants which
are intimately tied to the smooth structure of the manifold $Y$,
and the basic point is to bring  the special 
features of the Seiberg-Witten and 
Chern-Simons geometry
into play. Indeed, our
  new invariants are profoundly different from
the classical equivariant theory. The main 
analytic reason for
it  is  the fact that the spectrum of the Hessian of
the Chern-Simons functional is infinite in both positive and
negative directions. 

In the classical  situation  with e.g.
$G=S^1$,  the theory built on $ES^1 =S^{\infty}$ can be
approximated by the theories based on  
 $S^{2n+1}$. 
In spirit, our equivariant Seiberg-Witten Floer theories based
on $S^{2n+1}$ are also  approximations of 
 our theory based
on $S^{\infty}$.   
We shall refer to these equivariant theories 
(with $S^{2n+1}$ or $S^{\infty}$) as 
the {\it singular version} of equivariant Seiberg-Witten 
Floer homology and cohomology.

Instead of chains and cochains on the quotient space
$\R^0 \times \SS$, we can also use invariant cochains on
the stablized space $\R^0 \times \SS$. This is the
formulation presented in $\WangYeI$, an earlier 
version of [$\WangYe$],
which we call the 
{\it stable equivariant theory}.   It is essentially
equivalent to the above (quotient) equivariant
 construction. (At that time we already formulated 
 the quotient version, but didn't include it in the 
 paper. The reason was that we could only handle 
 $\SS=S^1$, and hence the advantage of the 
 quotient version could not be exploited, 
 see the relevant discussion below.)           
The precise relations between the equivariant and
stable equivariant theories, and between the $S^{2n+1}$ 
theory
and the $S^{\infty}$ theory,  
will be computed in Part III $\Ye$.  We acknowledge
a helpful  conversation with K. Fukaya on these
concepts.

\bigskip
\noindent {\bf Invariance}
\smallskip

The most fundamental issue about these theories
 is their diffeomorphism invariance. We 
establish their  invariance
by utilizing the {\it stablized transition Morse-Floer-Bott
flow complex} and its $S^1$ quotient, the
{\it quotient stablized transition Morse-Floer-Bott 
flow complex}. (In the proof, certain extensions 
of the transition flow allowing additional parameters
are needed. We  mean to include them here.) 
Note that  here stablization is not simply 
taking product with $\SS$. Instead, $\SS$ is built into 
the stable transition flow equation in a nontrivial 
fashion. As in [$\WangYe$], one has to kill
 reducible
transition trajectories in order to achieve
transversality.
 For this purpose, 
we employ    a  key  {\it gauge equivariant spinor
 perturbation}
for the stable transition flow equation. 
It is 
a modification of the   spinor perturbation 
for the transition flow equation 
introduced in  
[$\WangYe$], which is  equivariant only with 
respect to based gauges.   
 Note that for the Bott-type theory in
[$\WangYe$], the spinor perturbation of the 
transition flow equation provides only part of the 
invariance proof.  However, the gauge 
equivariant spinor perturbation of the stable 
transition flow equation suffices fully for 
establishing the invariance of the 
equivariant theories.    

We would like to mention that 
in  $\WangYeI$, besides the spinor perturbation
for the transition flow 
equation which is based gauge equivariant,
 we already constructed a 
gauge equivariant spinor perturbation for the 
stable transition flow equaion. However, it worked only 
for the case $\SS = S^1$. 
Now we have  a   
way to extend our construction to the cases 
$\SS= S^{2n+1}, n=1,2,...$ 
and 
$\SS=S^{\infty}$, and hence we are able to 
complete our 
equivariant theory. The completion of the singular version 
also allows us to establish the invariance of the de Rham 
version and the Cartan version, which will be addressed below.

There is another subtle point in the invariance 
proof we would like to mention. In the second step of 
the proof, we utilize a certain extension $IIE_{\e}$
of 
the transition flow equation ($\e$ is a parameter)
to establish the 
isomorphism equation
$F^+_{-*} \cdot F^-_{+*} =Id$. It turns out that
$IIE_0$ coincides with the original Seiberg-Witten flow 
equation.  The delicate feature here is that
all the moduli spaces $\M_T^0(S_{\a} \times_{S^1} \bold S, 
S_{\b} \times_{S^1} \bold S) \equiv \M^0_T(S_{a}, S_{\b})
\times_{S^1} \bold S$ with $\a =\b$ as well 
as $\a \not = \b$ are involved,
but no time translation is allowed. At a first glance,
it appears that we are running into trouble with 
compactification, for we had to use  the time translation
quotient $\underline \M_T^0(S_{\a}
\times_{S^1} \bold S, S_{\b} \times_{S^1} 
\bold S)$ in order to 
obtain the compactification $\underline \MM_T^0(S_{\a}
\times_{S^1} \bold S, 
S_{\b} \times_{S^1} \bold S)$. This trouble is resolved when 
we realize that
here a different compactification scheme is at play. 
As a consequence, the contribution of the moduli 
spaces $\M_T^0(S_{\a} \times_{S^1} S^1, S_{\b}
\times_{S^1} S^1)$ 
with $\a  \not = \b$ to 
the induced  chain map is trivial, and the 
remaining moduli spaces $\M_T^0(S_{\a} \times_{S^1} S^1,
S_{\a} \times_{S^1} S^1)$ give rise to the identity 
on the right hand side of  the equation 
$F^+_{-*} \cdot F^-_{+*} =Id$. 
Note that this 
 feature is irrelevant in Floer's situation $\FloerII$,
 where only 
 moduli spaces corresponding to $\a =\b$ are involved.

\bigskip
\noindent {\bf de Rham Version}
\smallskip

Instead of generalized cubical 
singular chains and cochains on the quotient 
space $\R^0 \times_{S^1} \SS$, we can also use
 differential forms on it.
The resulting theory will be called the {\it de Rham version of 
equivariant Seiberg-Witten Floer homology and cohomology}.
The construction of the de Rham  version
uses in  an essential way the fibration property of the 
projections of the quotient stablized  Morse-Floer-Bott
 flow complex. 
To establish
the invariance of this version, one might want to use 
the quotient stablized  transition Morse-Floer-Bott flow 
complex. However, one has  to
perform perturbations to kill reducilbe 
transition trajectories in order to bring this flow 
complex into general 
position.  The key perturbations  we
use for proving the invariance of the singular version achieve
these goals, but the resulting projections of the transition flow
complex may not be fibrations. 
In general, we think that it is unlikely to 
find suitable perturbations without 
destroying the fibration property of the 
projections.
Thus, it seems that any attempt of directly proving the 
invariance of the de Rham version in terms of 
the transition flow is doomed to fail.

We shall establish the invariance of the de Rham version
by showing   that it is isomorphic to the singular
version with real coefficients. The tool for establishing 
the isomorphism is the  spectral sequence induced 
from the  index filtration associated with the Chern-Simons
functional.

We shall also construct the de Rham version of the Bott-type 
Seiberg-Witten Floer homology and cohomology and show that 
it is isomorphic to the singular version constructed in
Part I [$\WangYe$] with real coefficients.   The relations
 between the equivariant 
theories and the Bott-type theories will be computed in
Part III.

\bigskip
\noindent {\bf  Cartan
Version} \smallskip

Besides Borel's homotopy quotient construction,
there is another  well-known  classical construction of 
equivariant cohomology based on equivariant  Lie algebra valued
differential 
forms.  It  is due to H. Cartan $\Cartan$ and  isomorphic to 
the Borel construction with real coefficients. 
In $\AustinBraamI$, equivariant instanton Floer homology 
and cohomology were constructed by utilizing 
Cartan's model.  Similar constructions 
can be carried out in the Seiberg-Witten set-up.
A number of  aspects of this construction have
independently  been 
presented in $\MarcolliWang$.  We present an account of 
this  construction, which is based on
our theory of the Morse-Floer-Bott flow complex along 
with its projection to critical submanifolds.  
This version will be referred to as the {\it Cartan 
version  of equivariant 
Seiberg-Witten Floer homology and
cohomology}.  A major difference between this version 
and
the singular  version as well as the  de Rham version of 
equivariant theory  is 
that here the chains and
cochains live on the moduli space $\R^0$ rather than the
stablized  moduli space $\R^0 \times S$ or its $S^1$
quotient.  

Like the de Rham version,  the construction of the 
Cartan  version
uses in  an essential way the fibration property of the 
projection of the {\it unstablized} Morse-Floer-Bott flow complex. 
Similar to the situation of the de Rham version,
one runs into troubles when trying  to prove its 
invariance directly by using
the transition Morse-Floer-Bott flow complex. Indeed, the
trouble here goes deeper than
 the fibration property. In order to 
deal with  equivariant differential forms, one 
needs to retain 
the full gauge equivariance of the transition flow 
equation  when perturbing it into transversal position.
Unfortunately, this is impossible to achieve. 
(Although we can construct full gauge 
equivariant transversal perturbations for 
the {\it stable} transition flow equation as discussed 
above.) Otherwise, one 
would be able to use the constructions to show that the 
conventional Serberg-Witten Floer homology is invariant, which,
according to $\DonaldsonII$, $\Chen$ and $\Lim$, fails to be 
true for homology spheres in general.  
(The paper $\MarcolliWang$  misses this crucial  point.
Indeed, in $\MarcolliWang$  
the issue of  transversal perturbation for the transition 
flow is not addressed.)

We  prove the invariance of the Cartan version
by  showing that it is isomorphic to the de Rham version.
   The isomorphism will be established
by  utilizing the  spectral
sequence induced from the index filtration associated
with the Chern-Simons functional.

\bigskip 
\noindent {\bf Exponential Convergence}
\smallskip

We present a detailed treatment of the 
issues of exponential decay and convergence 
of Seiberg-Witten trajectories. The convergence 
analysis contains  four themes.
 First, uniform 
pointwise 
smooth estimates hold for  a sequence of 
temporal Seiberg-Witten trajectories of
 uniformly 
bounded energy, provided that they are in 
appropriate gauges. Using this we obtain 
local convergence of the trajectories. 
But energy might be lost near infinities
during convergence,
which is not captured by the limit of local 
convergence. We manage to get additional limits
and show that they together capture all 
energies. This is the second theme and 
 quite similar in spirit to the convergence 
arguments about pseudo-holomorphic curves 
presented in $\YeGromov$. The third theme is 
to show that the convergence of the 
relevant piece of the trajectory 
to each limit actually holds in exponential 
norm.  This is an indispensible point.
 The fourth theme
is to show that the  endpoints of the limit 
temporal
trajectories match each other.
In other words, the endpoint at the positive time infinity
of one limit trajectory is identical to the endpoint of 
the next limit trajectory at the negative time 
infinity.  Indeed, this property  lies at the heart of 
the analysis of the structures of the compactified 
moduli spaces of trajectories.  

We would like 
to point out  that in the above analysis it is 
important to deal with temporal trajectories, i.e. 
use temporal gauges. 
The first reason is that  the key 
analytical properties of the trajectories
can be best utilized in temporal gauges.
(We also need to use local Columb gauges and 
``slice gauges", namely gauges determined 
along each $Y$-slice, as auxilliary tools.) 
The second reason is that the 
projections of our moduli spaces of 
trajectories have to be defined in terms of temporal 
gauges, see [$\WangYe$]. (Indeed they  are called 
{\it temporal projections}.) For this reason,
it is crucial to establish the 
endpoint matching property  
for limit temporal trajectories.   
Further relevant discussions are given 
below.

\bigskip 
\noindent {\bf Structure Near Infinity}
\smallskip

The convergence analysis discussed above is one pillar 
in the  analysis of the structures of the (compactified) 
moduli spaces of trajectories. The goal of this 
analysis is to show that these moduli spaces are 
compact smooth manifolds with corners.  The other pillar is 
a procedure of 
gluing trajectories, i.e. deforming piecewise trajectories 
into smooth trajectories.  Gluing arguments 
have  widely been used in gauge theories, in particular 
in Floer's instanton homology theory. There are delicate 
new aspects in our situation, however. First,
as mentioned above and explained in [$\WangYe$],
the endpoint projections have 
to be defined in terms of the temporal model.  
The appropriate objects as compactification
limits of trajctories are then the {\it consistent 
piecewise trajectories} and 
{\it consistent multiple 
temporal trajectory classes}, where ``consistant" means that 
the temporal projections of the pieces in a piecewise 
trajectory (or multiple 
trajectory class) match each other.  Our convergence 
argument 
shows that  trajectories in suitable gauges 
indeed converge to 
consistent piecewise trajectories.
Now we need to show the converse, namely consistent 
piecewise trajectories can be deformed back into 
smooth trajectories. This is done by a carefully 
designed gluing 
process. The crucial  new feature here 
is  the consistent 
condition. Indeed, conventional  gluing 
set-ups in the literature 
do not take care of the temporal  projections which are 
the key in the constructions of the equivariant and 
Bott-type theories. We carry out the gluing process 
in the set-up of the temporal model $\underline 
\MM_T^0(S_{\a}, S_{\b})$, which we 
prefer because of its canonical formulation 
 and global features.
 We also 
sketch this process  in the set-up of
 the fixed-end model $\underline \MM^0(p, q)$.

The second  delicate aspect is this.
Merely gluing piecewise trajectories into 
smooth trajectories falls far short from 
establishing the structure of smooth manifolds 
with corners. What is needed is suitable 
coordinate charts based on the gluing 
construction. To show that the gluing 
construction indeed yields the desired 
charts, we 
need to show 
the local diffeomorphism (in the interior) and  
homeomorphism property of 
the gluing construction. 
For this purpose, it is crucial to derive 
careful estimates for 
various sizes involved in the 
gluing construction. (There are 
treatments of this issue in the literature
with various degrees of details,
but our situation is very different.)

\head{2. Singular Versions of equivariant theory}
\endhead 

Recall the set-up in [$\WangYe$]: we consider a 3-dimensional 
rational homology sphere $Y$, along with a Riemannian metric $h$ and 
a $spin^c$ structure $c$ on $Y$.  
(For simplicity, we only 
present the case that $Y$ is connected. Our 
arguments work equally well for $Y$ with more than 
one components.)
Moreover, we consider  a pair of $Y$-generic 
parameters $(\pi, \lambda)$, namely
generic parameters  
for the perturbed Seiberg-Witten 
equation (1.1).
Recall also that for $\a \in 
\R$, $S_{\a}$ means the lift of $\a$ to 
$\R^0$, and that $S_i = \cup\{S_{\alpha}:
\alpha \in \R, \mu(\alpha)=i\}$.   
(The index $\mu$ was defined in  [$\WangYe$, Section 5].) 
 We have 
 $\R^0 = \cup_i S_i$.
% $S_i = \cup\{S_{\alpha}:
%\alpha \in \R, \mu(\alpha)=i$. 

Finally, we recall that (see [$\WangYe$, Section 7]) 
 to a topological space 
$\Cal X$ and coefficient group $\GG$, we associate 
the complex    $(C_*(\Cal X;\GG),
\pa_O)$ of generalized cubical singular 
chains, and the complex $
(C^*(\X; \GG),
\pa_O^*)$ of 
generalized cubical singular cochains.

\bigskip

\noindent {\bf Equivariant  Seiberg-Witten 
Floer Homology and Cohomology} 

\smallskip
 
Let $\SS$ be a topological 
space with a free $S^1$ action.  We choose $\SS$ to be 
either the odd dimensional euclidean spheres 
$S^{2n-1}   \subset \Bbb C^n$ or $S^{\infty}$.
The action of $S^1$ on them is multiplication by 
unit complex numbers.
Note that there are two well-known models for $S^{\infty}$. 
One is the unit 
sphere in a seperable Hilbert space. The other is the  
unit sphere in the space $\Bbb C^{\infty}$ whose elements 
are sequences $(c_1, c_2, ...)$ of complex numbers with finitely 
many nonzero  entries. The second model is precisely the 
direct limit of $S^{2n-1}$ as $n$ approaches infinity. 
We can use either one to carry out our constructions. 
We choose the second one, because it is more convenient. 

 We have the diagonal 
action of $S^1$ on $\R^0 \times \SS$ and the corresponding 
quotient 
$\R^0 \times_{S^1} \SS = \cup_i S_i \times_{S^1} \SS$. (Recall that the 
action of 
$S^1$ on $\R^0$ is via the identification of $S^1$ with constant 
gauges.)
Fix a coefficient group $\GG$. 
We introduce  our equivariant chain complex 
$C^{equ}_* =C^{equ; \GG}_*$
and cochain complex $C^*_{equ} = C^*_{equ;\GG}$,
$$C^{equ}_k=\oplus_{i+j=k}  C_j(S_i \times_{S^1} \SS; 
\GG),
C^k_{equ}=\oplus_{i+j=k} 
C^j(S_i \times_{S^1} \SS; \GG).$$
We have  $C^{equ}_* =  C_*(\R^0 \times_{S^1} \SS; \GG),
C^*_{equ;\GG}= C^*(\R^0 
\times_{S^1} \SS; \GG)$.

To construct the desired boundary operators for 
the equivariant chain and cochai-
n complexes, we consider 
the   moduli spaces of  consistent multiple temporal 
Seiberg-Witten trajectory classes $\underline { \MM}_T^0
(S_{\a}, S_{\b})$ introduced 
in [$\WangYe$, Section 6]  and their stablization 
$\underline {\MM}_T^0(S_{\a}, S_{\b})\times \SS$.
We write the latter as $\underline \MM_T^0(S_{\a} \times
\SS,  S_{\b} \times \SS)$. Indeed, a pair $(u, \ss)$ 
with $u$ a Seiberg-Witten trajectory can be considered 
as a solution of the trivially stablized 
Seiberg-Witten trajectory equation, namely the Seiberg-Witten
trajectory equation  
with a non-appearing, trivial variable  $\ss$.  With 
this set-up,
the moduli spaces  $\underline 
\MM_T^0(S_{\a}, S_{\b})$ 
used in [$\WangYe$] for the Bott-type theory 
are replaced 
by  $\underline \MM_T^0(S_{\a} 
\times \SS, S_{\b} \times \SS)$.  For this reason, 
we adopt the notation $[u, s]^T_0$ for 
$([u]_0^T, s)$.  (Recall that the subscript
$0$ refers to based gauges, while 
the  superscript $T$ means temporal.)

 We have the 
diagonal action of $S^1$ on these moduli spaces. 
(Recall that the action of 
$S^1$ on the Seiberg-Witten trajectories is via the identification 
of $S^1$ with constant gauges. ) Their quotients 
$\underline {\MM}_T^0(S_{\a} \times \SS, S_{\b}
\times \SS)/S^1$, which we shall denote by 
$\underline {\MM}_T^0(S_{\a} \times_{S^1} \SS, 
S_{\b} \times_{S^1} \SS)$,
now play the role of Morse-Floer-Bott flow complex for the 
equivariant 
chain complex and cochain complexes.  We call them the
 {\it quotient stablized  
Morse-Floer-Bott flow complex}. 
Since the temporal projections 
$\pi_+$ and $\pi_-$ in [$\WangYe$, Section 6] 
are equivariant with respect to the $S^1$ action,
we can define the  temporal  projections  for 
the quotient stablized 
 Morse-Floer-Bott flow complex  as follows
$$
\pi_+([u,\ss]_0^T)= ([\pi_+(u)]_0, \ss),  \pi_-([u, \ss]^T_0)=
 ([\pi_-(u)]_0,\ss).
$$
 With
these preparations, it is easy to carry over
 the construction of  
the boundary operators $\pa_{Bott}, \pa^*_{Bott}$ in
 [$\WangYe$, Section 7] 
to the present situation
to produce a  boundary operator $\partial_{equ}:
C^{equ}_k \to C^{equ}_{k-1}$ and the corresponding 
coboundary operator $\pa^*_{equ}: 
C_{equ}^{k} \to C_{equ}^{k+1}$.
Note that in the case $\SS= S^{\infty}$, the construction of 
these operators involves for each 
generalized cubical singular chain a fixed 
finite dimensional sphere
in $S^{\infty}$, and hence is no different from 
the construction in the case $\SS= S^{2n+1}$. This is because 
that every generalized singular cube 
has compact image in $\R^0 \times_{S^1} S^{\infty}$.
 (In general, for similar reasons,
all of our constructions
in the case  $\SS=S^{\infty}$ involve the same 
analysis as those in the case  
$S^{2n+1}$.)
Similar to $(\pa_{Bott})^2=0, 
(\pa^*_{Bott})^2=0,$ we have 
$(\partial_{equ})^2 =0, (\partial_{equ}^*)^2=0$.
Hence we can introduce the following equivariant 
Seiberg-Witten Floer homology
and cohomology.

\definition{Definition 2.1} We define (for each given 
coefficient group $\GG$ and $spin^c$ structure 
$c$) $$FH^{SW}_{equ*}
 =H_*(C^{equ}_*, \pa_{equ}), 
FH^{SW*}_{equ} = H^*(C^{equ}_*, \pa_{equ}^*).$$ 
\enddefinition

If we need to indicate  the coefficient group $\GG$
and the $spin^c$ structure $c$,, 
we can write e.g. $FH^{SW}_{equ; \GG*}(c)$.

As in [$\WangYe$], we  consider the index filtration 
$$
\F_*^{equ}=\cdot \cdot
\cdot \F_k^{equ} \subset \F_{k+1}^{equ}
\cdot \cdot  \cdot 
$$
 and its dual filtration 
$\F^*_{equ}$ for 
our equivariant chain and cochain complexes , where 
$$
\F_k^{equ} = \oplus_{j \leq k}  C_*(S_j \times_{S^1} \SS;
\GG),
$$
$$
\F^k_{equ} = \oplus_{j \geq k}  C^*(S_j
\times_{S^1} \SS; \GG).
$$

We have

\proclaim{Theorem 2.2} The equivariant index filtration
 $\F_*^{equ}$
induces a spectral 
sequence $E(equ)^*_{**}$ converging to $FH^{SW}_{equ*}$  such that
$$
E(equ)^1_{ij} \cong H_{j}(S_i \times_{S^1} \SS; \GG). 
\tag 2.1
$$
The dual filtration induces a spectral sequence 
$E(equ, dual)^*_{**}$ converging to  $FH^{SW*}_{equ}$
such that 
$$
E(equ, dual)^1_{ij} \cong H^j(S_i \times_{S^1} \SS; \GG).
\tag 2.2
$$
\endproclaim

The proof is similar to  the proof for  
[$\WangYe$, Theorem 7.8], hence we omit it.

\bigskip

\noindent {\bf Stable Equivariant Seiberg-Witten
Floer Homology and Cohomology}
\smallskip
We consider the same $\SS$  and fix a coefficient group 
$\GG$ as before.
The action of $S^1$ on $\R^0 \times \SS$ induces an action
on
generalized cubical singular chains in $\R^0 \times \SS$:
for generators $(\Delta, f)\otimes \gg$ and $g 
\in S^1$,
we have $g^*((\Delta, f)\otimes \gg) 
=(\Delta, g^*f(\cdot))\otimes \gg$.
Passing to quotients, we
obtain an action of $S^1$ on
$ C_*(\R^0 \times \SS)$.  A
$j$-cochain class
$\omega  \in  C^*(\R^0
\times \SS; \GG)$ is called {\it invariant},
provided that $\omega(g^*\sigma)=\omega(\sigma)$ for
all $\sigma \in C_*(\R^0 \times \SS; \Bbb Z)$ 
and $g \in S^1$.

We define $ C^j_{inv} (\R^0 \times \SS)$
to be the group of
invariant  $j$-cochain classes on
$\R^0\times S$ and introduce the stable equivariant 
cochain complex $C^*_{sequ}$, 
$$C^k_{sequ}=\oplus_{i+j=k} C^j_{inv}(S_i \times
 \SS).$$

Using the stablized moduli spaces $\underline 
{\MM}_T^0(S_{\a}, S_{\b}) \times \SS$
we can easily carry over the construction of 
the coboundary operator  
$\pa^*_{Bott}$  to obtain 
a boundary operator $\pa_{sequ}^*: C^k_{sequ}
\to C^{k+1}_{sequ}$ for the stable 
equivariant complex. 
The gauge equivariance of the projections $\pi_{\pm}$ 
ensures that $\pa_{sequ}$ indeed produces invariant 
cochains out of invariant cochains.
We have $\pa_{sequ}^2=0$.

\definition{Definition 2.3} The stable
equivariant Seiberg-Witten Floer
cohomology
$FH^{SW*}_{sequ}$
(for the  given coefficient group
$\GG$ and $spin^c$ structure $c$) 
is defined to be the 
cohomology 
of the cochain complex $(C^*_{sequ},
\pa_{sequ}^*)$. 
The stable equivariant Seiberg-Witten
Floer homology
$FH^{SW}_{sequ*}$ is defined
to be its dual homology. 
\enddefinition 

We have a stable equivariant   index filtration
$$\F^*(sequ)=\cdot \cdot \cdot \F^k(sequ) 
\subset \F^{k-1}(sequ) \subset \cdot    \cdot   
\cdot,$$
$$\F^k(sequ) = \oplus_{j \geq k} C^*_{inv}
(S_k \times  \bold S).$$
There is a corresponding  dual filtration. The 
corresponding spectral sequences are similar to 
those in Theorem 2.3.

\definition{ Remark 2.4} We can also formulate a stable version of the 
Bott-type theory in [$\WangYe$], whose information 
is contained in  the Bott-type theory. We leave the 
details to the reader.
\enddefinition

\head{3. Invariance }
\endhead

The purpose of this section is to prove 
the following theorem. 
        
\proclaim{Theorem 3.1} The singular version of 
equivariant  
Seiberg-Witten Floer homology
and cohomology are diffeomorphism 
invariants up to shifting isomorphisms.
\endproclaim

We start with 

\definition{Definition 3.2} 
Consider  $\SS=S^{2n-1}$ or $S^{\infty}$.
Let the group of based gauges $\G_{3, loc}^0 
\equiv \G_{3, loc}^0(X)$ (recall 
$X = Y \times \Bbb R$) act on 
$(\A_{2, loc}(X)
\times \Ga^+_{2, loc}(X)) \times \SS$
in the following fashion:
$$
g^*(u, \ss) = (g^*u,
\ss),
$$
where   $g\in
\G_{3, loc}^0,  u \in \A_{2, loc} 
\times \Ga^+_{2, loc}$ and $ \ss \in \SS$.  On the other hand,
we have the diagonal action of $S^1$ on 
$(\A_{2, loc} \times \Ga^+_{2, loc}) \times \SS$:
$$
g^*(u, \ss) = (g^*u, g^{-1} \ss),
$$
where $g \in S^1$. Combining these two actions we 
then obtain an action of the full group of gauges 
$\G_{3, loc}$ on $(\A_{2, loc} \times \Ga^+_{2, loc})
\times \SS.$              
\enddefinition

We first consider the case $\SS=S^1$, which allows 
a special treatment, simpler than the general one.

\definition{Definition 3.3} We define a smooth vector field
$Z_e$ on $(\A_{2, loc} \times\Ga^+_{2, loc}) \times S^1$
as follows $$Z_e(u, s) = s Z(u),  $$
where $Z$ is the vector field given by [$\WangYe$,
Definition 8.3].   
\enddefinition

The following lemma is readily
proved.
\proclaim{Lemma 3.4} $Z_e$ is equivariant with respect to
the action of $\G_{3, loc}$.
\endproclaim

For two metrics  $h_{\pm}$ on $Y$ and two 
pairs of $Y$-generic parameters $ 
(\pi_{\pm}, \lambda_{\pm})$ for $h_{\pm}$ 
respectively, we have 
their  interpolations $h(t), \pi(t),
\lambda(t)$ as given in [$\WangYe$, Section 8].
We introduce the following 
stable version of the (perturbed)   
 transition trajectory equation    
([$\WangYe$, (8.1)]) 
for $A= a + f dt, \P=\p
$ and $ s \in S^1$.  Its solutions will be called {\it stable
 transition 
trajectories}.

$$
\cases
\frac{\partial a}{\partial t} & =*_YF_a+d_Yf+
\langle e_i\cdot\p,\p\rangle e^i+
\n H_{\pi(t)}(a) +  b_0,
\cr
\frac{\partial \p}{\partial t} & =-\np_a\p-\lmd(t)\p+
 Z_{e}((A, \P), s),\cr
\endcases
\tag 3.1
$$
where $s \in S^1$,  and $ b_0$ is a 
perturbation form. Several operations in this 
equation (such as 
the Hodge star)
depend in an obvious way 
on the interpolating metric $h(t)$ at time $t$.
Note that this equation can be written in the 
following fashion:
$$
\cases
F^+_A =&\frac{1}{4}\langle e_ie_j\P,\P\rangle e^i\wedge e^j
\cr 
&+ ( \nabla H_{\pi(t)}
(a) +b_0) \wedge dt+*(  (\nabla H_{\pi(t)} (a)+
b_0) \wedge dt),\cr
D_A\P =& - \lambda \frac{\pa}{\pa t} \cdot 
(\P- Z_e((A, \P),s),\cr
\endcases
\tag 3.2
$$
where $X$ is endowed with the warped-product 
metric determined by $h(t)$ and the standard 
metric on $\Bbb R$. 

As a consequence of Lemma 3.4 and the construction of 
$Z_e$ we have
\proclaim{Lemma 3.5} The equation (3.1) is invariant under 
the action of $\G _{3, loc}$. Moreover, it has no 
reducible solution.
\endproclaim

We have moduli spaces of temporal stable 
transition trajectories $ 
\M_T^0(S_{\a_-}\times S^1, 
S_{\a_+} \times S^1)$ with $\a_- \in \R_-$
and $\a_+ \in \R_+$, where $\R_-$ and $\R_+$ are the 
moduli spaces of gauge equivalence classes of 
Seiberg-Witten points for the parameters
$(h_-, \pi_-, \lambda_-)$ and $(h_+, \p_+,
\lambda_+)$ respectively. They  are analogous 
to the moduli spaces $\M_T^0(S_{\a_-},
S_{\a_+})$ in [$\WangYe$, Section 8]. 
We also have the 
moduli spaces of stable, consistent  multiple
temporal transition trajectory classes (SCMTC), $\MM_T^0(
S_{\a_-} \times S^1, S_{\a_+} \times S^1)$. 
A SCMTC is an element $([u_1, s_1]_0^T,...,[u_k,
s_k]_0^T)$ 
in $\M^0_T(S_{\a^-_0}\times S^1, S_{\a_1^-}\times S^1)
\times_{S_{\a_1^-}}... \M^0_T(S_{\a_{m-1}^-}\times S^1,
S_{\a_m^+}\times S^1) \times_{S_{\a_m^+}}...
\M^0_T(S_{\a^+_{k-1}} \times S^1, S_{\a^+_k}
\times S^1)$. The time translation
acts on all portions except the distinguished 
one $[u_m, s_m]^T_0$.  The quotient of $\MM_T^0(S_{\a_-}
\times S^1, S_{\a_+} 
\times S^1)$ under the 
resulting $\Bbb R^k$-action is denoted by 
$\underline \MM_T^0(S_{\a_-}
\times_{S^1}  S^1, S_{\a_+} \times_{S^1} S^1)$. All these 
are similar 
to the constructions and notations in [$\WangYe$, Section 8
and Section 6].

We have the following analogue of [$\WangYe$, 
Theorem 8.9].

\proclaim{Theorem 3.6} The moduli spaces 
${\underline \MM}^0_T(S_{\a_-}
\times S^1,
S_{\a_+}\times S^1)$ are compact 
Hausdorff spaces with respect to 
the topology 
induced from smooth
piecewise exponential convergence 
(cf. [$\WangYe$, Definition 6.14]).
Moreover, for generic $ b_0$, the following
hold 
for all $\a_- \in \R_-, \a_+ \in \R_+$.

(1) ${\underline {\MM}}_T^0(S_{\a_-} \times S^1,
S_{\a_+}\times S^1)$ has the structure of
$d$-dimensional smooth oriented manifolds with
corners,
where $d = \mu_{-}(\a_-)- \mu_+(\a_+) +m_0 
+\hbox{dim }G_{\a_+}-max\{\hbox{dim }
G_{\a_-}, \hbox{dim }G_{\a_+}\} +2$, and 
$$m_0 
={ind}~ \F_{O_-, O_+}-1,$$ with $\F_{O_-, O_+}$
denoting the index of the linearized transition 
Seiberg-Witten operator between the two reducible 
Seiberg-Witten points $O_-, O_+$.
(The stable parameter in $S^1$ does not enter 
into the 
index counting. But it contributes an additional 
one to the dimension formula.)

(2) This structure
is compatible with the natural stratification of
$\underline {\MM}_T^0(S_{\a_-} \times S^1,
S_{\a_+}\times S^1)$.

(3) The (canonically defined) temporal 
projections $\pi_{\pm}:
{\underline {\MM}}_T^0(S_{\a_-} \times S^1,
S_{\a_+}\times S^1) \to  S_{\a_{\pm}} \times S^1$
are $S^1$-equivariant smooth maps.
 %where the action of $S^1$ on the 
%moduli spaces is induced by the following action: $s \in S^1$ 
%acts on $(u, s')$ to yield $(s^*u, s^{-1}s')$.
(But they may not be  fibrations in general.)

Consequently, there holds
$$
\pa \underline \MM_T^0(S_{\a_-}\times S^1, S_{\a_+}
\times S^1) = 
(\cup_{\mu(\a_-) > \mu(\a_-') \geq \mu(\a_+)-m_0}
\underline \MM_T^0(S_{\a_-} \times S^1, S_{\a_-'}
 \times
S^1)
\tag 3.3
$$
$$ 
\times_{(S_{\a_-'} \times S^1)} \underline \MM_T^0(S_{\a_-'}
\times S^1, 
S_{\a_+} \times S^1) )\cup (\cup_{\mu(\a_-) \geq 
\mu(\a_+')-m_0 > \mu(\a_+)-m_0} \underline \MM_T(S_{\a_-}
\times S^1, $$
$$
S_{\a_+'}\times S^1) \times_{(S_{\a_+'}\times 
S^1)} \underline \MM_T^0(S_{\a_+'} \times S^1,
S_{\a_+} \times S^1)).
$$
\endproclaim

Passing to $S^1$-quotient, we 
obtain  the following result.

\proclaim{Theorem 3.7} 
The moduli spaces 
${\underline \MM}^0_T(S_{\a_-}
\times_{S^1} S^1,
S_{\a_+}\times_{S^1} S^1)$ are compact
Hausdorff spaces. 
Moreover, for generic $ b_0$, the following
hold 
for all $\a_- \in \R_-, \a_+ \in \R_+$.

(1) ${\underline {\MM}}_T^0(S_{\a_-} \times_{S^1} S^1,
S_{\a_+}\times_{S^1} S^1)$ has the structure of
$d$-dimensional smooth oriented manifolds with
corners,
where $d = \mu_{-}(\a_-)- \mu_+(\a_+) +m_0 
+\hbox{dim }G_{\a_+}- max\{\hbox{ dim }
G_{\a_-}, \hbox{dim }G_{\a_+} \} +1$.

(2) This structure
is compatible with the natural stratification of
$\underline {\MM}_T^0(S_{\a_-} \times_{S^1} S^1,
S_{\a_+}\times_{S^1} S^1)$.

(3) The (canonically defined) temporal
projections $\pi_{\pm}:
{\underline {\MM}}_T^0(S_{\a_-} \times_{S^1} S^1,
S_{\a_+}\times_{S^1} S^1) \to  S_{\a_{\pm}} \times_{S^1} S^1$
are  smooth maps.
 %where the action of $S^1$ on the 
%moduli spaces is induced by the following action: $s \in S^1$ 
%acts on $(u, s')$ to yield $(s^*u, s^{-1}s')$.
(But they may not be  fibrations in general.)

Consequently, there holds
$$
\pa \underline \MM_T^0(S_{\a_-} 
\times_{S^1} S^1, S_{\a_+}
\times_{S^1} S^1) = 
(\cup_{\mu(\a_-) > \mu(\a_-') \geq \mu(\a_+)-m_0}
\underline \MM_T^0(S_{\a_-} \times_{S^1} S^1, 
\tag 3.4
$$
$$S_{\a_-'}
\times_{S^1} S^1) 
\times_{(S_{\a_-'} \times_{S^1} S^1)} \underline
 \MM_T^0(S_{\a_-'} \times_{S^1} S^1, 
S_{\a_+} \times_{S^1} S^1) )\cup 
$$
$$(\cup_{\mu(\a_-) \geq 
\mu(\a_+')-m_0 > \mu(\a_+)-m_0} \underline \MM_T(S_{\a_-}
\times_{S^1} S^1,
S_{\a_+'} \times_{S^1} S^1) \times_{(S_{\a_+'}
\times_{S^1} S^1)}$$
$$ \underline \MM_T^0(S_{\a_+'} 
\times_{S^1} S^1,
S_{\a_+} \times_{S^1} S^1)).
$$
\endproclaim

\demo{Proof of Theorem 3.1 for $\SS=S^1$}
Consider the above set-up of two sets of 
parameters $(h_{\pm}, \pi_{\pm}, \lambda_{\pm})$.
Our goal is to show that the singular version 
of equivariant Seiberg-Witten Floer homology 
(cohomology) constructed in terms of $(h_-,
\pi_-, \lambda_-)$ is isomorphic to 
that constructed in terms of $(h_+, \pi_+,
\lambda_+)$. We present the case 
of homology, while the case of cohomology 
can be handled by a similar argument.

Consider the equivariant 
complexes $C^{equ-}_*$ and $C^{equ+}_*$ 
associated with the parameters $(h_-, \pi_-,
\lambda_-)$ and $(h_+, \pi_+, \lambda_+)$ 
respectively.  Let $ \bold F^-_+$ denote the collection 
of all projections $\pi_-: \underline \MM_T^0(S_{\a_-}
\times_{S^1} S^1, 
S_{\a_+} \times_{S^1} S^1) \to S_{\a_-} \times_{S^1}
S^1$ and $\pi_-:
\pa \underline \MM_T^0(S_{\a_-} \times_{S^1} S^1, S_{\a_+}
\times_{S^1} S^1) \to
S_{\a_-} \times_{S^1} S^1$ (for all $\a_-, \a_+$).
 We have  the subcomplex 
$C_*^{equ-, \bold F^-_+}$ of $C_*^{equ-}$ 
consisting of $\bold F^-_+$-transversal 
chains on $\R^0_- \times_{S^1} S^1$, which are 
similar to $C_*^{Bott,
\bold F}$ in [$\WangYe$, Section 7]. By the 
arguments in the proof of [$\WangYe$, Lemma 7.8],
the homology $H_*(C^{equ-, \bold F^-_+}, \pa_{equ})$
is canonically isomorphic to $H_*(C^{equ-}_*,
\pa_{equ})$.
                   
Using the moduli spaces 
$\underline \MM_T^0(S_{\a_-} \times_{S^1} S^1,
S_{\a_+} \times_{S^1} S^1)$,
we construct a chain map $F^-_+: C^{equ-}_* \to 
C^{equ+}_*$ in the same way as the construction of 
the chain map $F$ in [$\WangYe$, Section 8]. 
We denote the interpolations and perturbation parameters
involved in (3.1) by $P^-_+=(h^-_+, \pi^-_+,
\lambda^-_+,  b^-_{0+})$.
Reversing the roles of $(h_-, \pi_-, 
\lambda_-)$ and $(h_+, \pi_+, \lambda_+)$,
we obtain a chain map $F^+_-: C^{equ+, \bold F^+_-}_*
\to C^{equ-}_*$. The corresponding 
interpolations and perturbation parameters 
are denoted by $P^+_-=
(h^+_-, \pi^+_-, \lambda^+_-, 
b^+_-)$.  The associated collection of 
projections is denoted by $\bold F^+_-$. 
Furthermore, let $\bold F^-_+ 
\bold F^+_-$ denote the collection of 
projections $\pi_-$ from $\underline 
\MM_T^0(S_{\a_-} \times_{S^1} S^1, 
S_{\a_+} \times_{S^1} S^1) \times_{(S_{\a_+} \times_{S^1}
S^1)} \underline \MM_T^0(S_{\a_+} \times_{S^1} S^1,
S_{\a_-'} \times_{S^1} S^1)$ to $S_{\a_-} \times_{S^1}
S^1$.

 Note that 
$F^-_+$ has degree $m_0$, and $F^+_-$ has degree 
$-m_0$. 
We need to show that they induce isomorphisms 
between the homologies.  The arguments consist of 
three steps. 
\bigskip

\noindent {\bf Step 1} 

Consider $\tau_{-R^{-1}}(P^-_+)$ and 
$\tau_{R^{-1}}(P^+_-)$, the time translated 
$P^-_+$ and $P^+_-$ for  $R \in 
(0, 1)$. Note that 
they coincide over $Y \times [1-R^{-1}, -1+R^{-1}]$.
We define $P_-= P_-(R)$ to be 
equal to $\tau_{-R^{-1}}(P^-_+)$ over 
$Y \times (-\infty, -1+R^{-1}]$ and 
equal to $\tau_{R^{-1}}(P^+_-)$ over $Y \times 
[1-R^{-1}, \infty)$.
Using $P_-$ in (3.1) we obtain the 
{\it interpolated stablized transition 
equation}. For a fixed $R$, we denote it 
by $IE_R$. If we allow $R$ to vary in an 
interval $(0, R_0]$, we denote it by
$IE_{(0, R_0]}$.  Roughly speaking, 
when $R \to 0$, $IE_R$ breaks into 
the stablized transition equation with 
data
$P^-_+$ and the same equation with 
data $P^+_-$.

Fix some $R_0 \in (0, 1)$. Using $IE_{R_0}$ we 
obtain moduli spaces $\underline \MM_T^0(S_{\a_-}
\times_{S^1} S^1, S_{\a_-'} \times_{S^1} S^1)^{IE_{R_0}}$
along with the collection $\bold F_{IE_{R_0}}$ of  projections
$\pi_-$
from them and their boundaries.  In order 
to achieve transversality for these moduli spaces,
we add an additional transversal perturbation 
form $\bar b_0$ to $P_-(R_0)$.  Then a structure 
theorem similar to Theorem 3.7 holds for them.
Using these moduli spaces we construct a chain map 
$F_-: C^{equ-, \bold F_-}_*
\to C^{equ-}_*$ of degree zero.

On the other hand, using $IE_{(0, R_0]}$
we obtain moduli spaces $\underline 
\MM_T^0(S_{\a_-}
\times_{S^1} S^1, S_{\a_-'} $
$\times_{S^1} S^1)^{IE}$. 
We denote the associated collection of 
projections by $\bold F_{IE}$. To achieve
transversality for these moduli spaces, we add a transversal 
perturbation famility of  forms $ b_0(R) $ to  $P_-(R)$,
such that $b_0(R_0)=\bar b_0$ and $b_0(R) \to 0$ 
as $R \to 0$. Then a structure theorem similar 
to Theorem 3.7 holds for them.    
Employing these moduli spaces  we  obtain 
a chain map $\Theta:  C^{equ-, \bold
 F_{IE}}_* \to  C^{equ-}_*$ of 
degree one.

We have 
$$\pa \underline \MM_T^0(S_{\a_-} \times_{S^1} S^1, 
S_{\a_-'} \times_{S^1} S^1)^{IE} =
 (\cup_{\a_-''} \underline \MM_T^0(S_{\a_-} \times_{S^1} S^1,
S_{\a_-''} \times_{S^1} S^1 ) 
\times_{(S_{\a_-''} \times_{S^1}
S^1)}
\tag 3.5
$$
$$
 \underline \MM_T^0(S_{\a_-''} \times_{S^1}
S^1, S_{\a_-'} \times_{S^1} S^1)^{IE} ) \cup
$$
$$
(\cup_{\a_-''} \underline \MM_T^0(S_{\a_-}
\times_{S^1} S^1, S_{\a_-''} \times_{S^1} S^1)^{IE} 
\times_{(S_{\a_-''} \times_{S^1} S^1)} 
\underline \MM_T^0(S_{\a_-''} \times_{S^1} 
S^1, S_{\a_-'} \times_{S^1} S^1) ) \cup
$$
$$
\underline \MM_T^0(S_{\a_-} \times_{S^1} S^1, 
S_{\a_-'} \times_{S^1} S^1)^{IE_{R_0}} \cup
$$
$$
(\cup_{\a_-''} \underline \MM_T^0(S_{\a_-} \times_{S^1} S^1,
S_{\a_-''} \times_{S^1} S^1) \times_{(S_{\a_-''} \times_{S^1}
S^1)} \underline \MM_T^0(S_{\a_-''} \times_{S^1} S^1,
S_{\a_-'} \times_{S^1} S^1)).
$$
Using this formula along with suitable orientations 
we can argue as 
 in the proofs of [$\WangYe$, Lemma 7.5 ]
and [$\WangYe$,  Theorem 8.10] to infer that 
$\Theta$  is a chain homotopy between 
$F^+_- \cdot F^-_+$ and $F_-$, i.e. $F_+ \cdot
F_- - F^-_+ =  \pa_{equ} \cdot \Theta + \Theta \cdot
\pa_{equ}. $ Here, we restrict the chain maps 
$F^-_+, F^+_-, F_- $ and $\Theta$ to 
the subcomplex $C^{equ, \bold F}_*$, where 
$\bold F=\bold F^-_+ \cup \bold F_-  \cup \bold F_{IE}
\cup \bold F^+_- \bold F^-_+$. 
It follows that for the induced maps on 
homologies:
$$
F^-_{+*}\cdot F^+_{-*}= F_{-*}.
$$

\bigskip

\noindent {\bf Step 2}

In the equation $IE_{R_0}$, we have a term 
corresponging to $b_0$ in (3.1), and a term 
$Z_e$. We multiply these two terms by a real parameter
$\e$ and obtain a new equation $IIE_{\e}$. Clearly,
$IIE_{1}=IE_{R_0}$. On the other hand, 
$IIE_0$  
is  the original perturbed Seiberg-Witten trajectory equation
(1.2) 
with some specific holonomy 
and $\lmd$ perturbations. 

We use $IIE_{[0, 1]}$ to denote the equation
$IIE_{\e}$ in which $\e$ is allowed to vary in 
$[0, 1]$.  Using it we obtain moduli spaces 
$\underline \MM_T^0(S_{\a_-} \times_{S^1}
S^1, S_{\a_-'} \times_{S^1} S^1)^{IIE}$ 
along with the associated collection of 
projections $\bold F_{IIE}$. On the hand,
we use $IIE_0$ to obtain moduli spaces 
$\underline \MM_T^0(S_{\a_-} \times_{S^1} 
S^1, S_{\a_-'} \times_{S^1} S^1)^{IIE_0}$.

The latter moduli spaces are not the same as 
$\underline \MM_T^0(S_{\a_-} \times_{S^1} S^1,
S_{\a_-'} \times_{S^1} S^1)$ used for the 
construction of our homology theory, althouth 
they are both built in terms of the solutions 
of the same equation (1.2).
{\it This distinction is essential for our purpose.} 
The reason for the distinction is as follows.
The moduli space $\M_T^0(S_{\a_-}
\times_{S^1} S^1, S_{\a_-'} \times_{S^1})^{IIE_0}$
is of course the same as $\M_T^0(S_{\a_-}
\times_{S^1} S^1, S_{\a_-'} \times_{S^1} 
S^1)$ for the same perturbation paremeters.
However, since $IIE_{\e}$ is not time translation 
invariant for $\e >0$, the 
$\Bbb R$-action on $IIE_{[0, 1]}$-trajectories 
is defined to be trivial. This dictates that 
we also define the $\Bbb R$-action on 
$IIE_0$-trajectories to be trivial. Hence the 
$\Bbb R$ action on $\M_T^0(S_{\a_-}
\times_{S^1} S^1, S_{\a_-'} \times_{S^1} S^1)^{IIE_0}$ 
is defined to be trivial. The compactification of 
this moduli space leads to $\underline
\MM_T^0(S_{\a_-}
\times_{S^1} S^1, S_{\a_-'} \times_{S^1} S^1)^{IIE_0}$.
We identify $\M_T^0(S_{\a_-} \times_{S^1} S^1,
S_{\a_-'} \times_{S^1} S^1)^{IIE_0}$ with $\underline 
\M_T^0(S_{\a_-} \times_{S^1} S^1, S_{\a_-'}
\times_{S^1} S^1) \times (-1, 1)$ via the 
time translation and an identification of 
$\Bbb R$ with $(-1, 1)$. Then we 
We have for 
$\a_- \not = \a_-'$
$$
\underline \MM^0_T(S_{\a_-} \times_{S^1} 
S^1, S_{\a_-'} \times_{S^1} S^1)^{IIE_0}=
\underline \MM^0_T(S_{\a_-}
\times_{S^1} S^1, S_{\a_-'} \times_{S^1} 
S^1) \times [-1, 1].
\tag 3.6
$$
Note that the projections $\pi_+$ on 
them are constant along the fiber $[-1, 1]$.

On the other hand, we have for $\a_-' =\a_-$
$$
\underline \MM^0_T(S_{\a_-} 
\times_{S^1} S^1, S_{\a_-} \times_{S^1}
S^1)^{IIE_0} = \M^0_T(S_{\a_-} \times_{S^1}
S^1, S_{\a_-} \times_{S^1} S^1).
\tag 3.7
$$
Indeed, here only time-independent 
trajectories are involved.

The transversality for 
$\M^0_T(S_{\a_-} \times_{S^1} S^1,
S_{\a_-} \times_{S^1} S^1)$ follows from [$\WangYe$,
Lemma C.1]. We incoporate  an additional holonomy
 perturbation to 
achieve transversality for $\underline \MM_T^0(S_{\a_-}
\times_{S^1} S^1, S_{\a_-'} \times_{S^1} S^1)$ for 
$\a_- \not = \a_-'$, cf. [$\WangYe$, Section 5].
Then we can incoporate additional holonomy and form 
perturbations to achieve transversality for 
$\underline \MM^0_T(S_{\a_-}
\times_{S^1} S^1, S_{\a_-'} \times_{S^1} S^1)^{IIE}$.
We have structure theorems for $\underline 
\MM_T^0(S_{\a_-} \times_{S^1} S^1,
S_{\a_-'} \times_{S^1} S^1)^{IIE_0}$ and 
$\underline \MM_T^0(S_{\a_-} \times_{S^1} 
S^1, S_{\a_-'} \times_{S^1} S^1)^{IIE}$ which 
are similar to Theorem 3.7.  We have 
the associated collections of projections 
$\bold F_{IIE_0}$ and $\bold F_{IIE}$.
We set $\bold F= \bold F_{IE_{R_0}} 
\cup \bold F_{IIE_0} \cup \bold F_{IIE}$.

Using the moduli spaces $\underline \MM_T^0(S_{\a_-}
\times_{S^1} S^1, S_{\a_-'} \times_{S^1} S^1)^{IIE_0}$
we construct a chain map $ F_0:
C_*^{equ-, \bold F} \to C_*^{equ-}$.  It is the sum of 
two maps, with one corresponding to $\a_- = \a_-'$,
and 
one corresponding to $\a_- \not = \a_-'$. By 
(3.6) and the fact that 
the projections $\pi_+$ are 
constant along the fiber $[-1, 1]$,
we deduce that
the image of the second map is given 
in terms of degenerate generalized singular cubes.
Hence this map is the zero map.
On the other hand, it follows from 
(3.7) that the first map is the inclusion 
map  ( the identity map). Hence $F_0$ is 
the inclusion map.

Using the moduli spaces $\underline 
\MM_T^0(S_{\a_-} \times_{S^1} S^1,
S_{\a_-'} \times_{S^1} S^1)^{IIE}$ 
we construct a chain map $\ti \Theta:
C_*^{equ-. \bold F} \to C_*^{equ-}$.
There holds
$$
\pa \underline \MM_T^0(S_{\a_-} \times_{S^1} S^1,
S_{\a_-'} \times_{S^1} S^1)^{IIE} = 
\tag 3.8
$$
$$
\underline \MM^0_T(S_{\a_-} \times_{S^1} 
S^1, S_{\a_-'} \times_{S^1} S^1)^{IIE_0} \cup
\underline \MM_T^0(S_{\a_-} \times_{S^1} S^1,
S_{\a_-'} \times_{S^1} S^1)^- \cup $$
$$
(\cup_{\a_-''} \underline \MM_T^0(S_{\a_-}
\times_{S^1} S^1, S_{\a_-''} \times_{S^1} S^1) 
\times_{S_{\a_-''} \times_{S^1} S^1} 
\underline \MM_T^0(S_{\a_-''} \times_{S^1} S^1,
S_{\a_-'} \times_{S^1} S^1)^{IIE}) \cup $$
$$
(\cup_{\a_-''} \underline \MM_T^0(S_{\a_-}
\times_{S^1} S^1, S_{\a_-''} \times_{S^1}
S^1)^{IIE} \times_{S_{\a_-''} \times_{S^1} S^1}
\underline \MM_T^0(S_{\a_-''} \times_{S^1} S^1,
S_{\a_-'} \times_{S^1} S^1).
$$

Using this equaton and suitable orientations we 
then deduce that $\ti \Theta$ is a chain homotopy between
$F_0$ and $F_-$, where $F_-$ is restricted 
to $C_*^{equ-, \bold F}$. It follows that $F_{-*}=Id$. 
 Hence $F^+_{-*} 
\cdot F^-_{+*} = Id$.

\bigskip

\noindent {\bf Step 3} 

A  similar construction with the roles of $ F_-$ and $F_+$ 
reversed yields $F^-_{+*} \cdot F^+_{-*}=Id$.
\qed
\enddemo

\bigskip

We have proved Theorem 3.1 for $\SS=S^1$. Now we handle 
the general case. For higher dimensional spheres, it's 
no longer possible to find an equivariant $S^1$ valued 
function, hence the spinor construction given in Definition 
3.2 cannot be applied. We introduce a  
new device.  Choose Hermitian orthonormal smooth 
spinor fields  
$\P_k, k=1, 2, ... \in \Ga^-(Y
\times [-1, 1]) \equiv \Ga(W^-|_{Y \times [-1, 1]})$ 
with supports contained in $Y \times (-1, 1)$. (See
[$\WangYe$] for the definition of the 
spinor bundles $W^+$ and $W^-$. )
Let $e_k$ be the vector in $\Bbb C^{\infty}$ whose 
$k$-th entry is 1, and all other entries are 
zero. The assignment $e_k \to \P_k$ determines 
an Hermitian complex linear 
embedding $\bold \Psi$ from $\Bbb C^{\infty}$ 
into $\Ga^-(Y \times [-1, 1])$. In particular, 
it is equivariant under the $S^1$ action,
which is multiplication by unit complex numbers.

\definition{Definition 3.8 } Consider the global slice 
$\S_1$ given in [$\WangYe$, Lemma 
8.1] for the action of $\G^0_3(X_1)$ on $\A_2(X_1)
\times \Ga^+_2(X_1)$. (Recall $X_1 = Y 
\times [-1, 1]$.) We define a smooth 
vector field $\ti Z$ on $(\A_2(X_1) \times \Ga^+_2(X_1))
\times \SS$ by
$$
\ti Z(g^*u, \ss) = g^{-1} \bold \Psi (\ss),
$$
for $g \in \G^0_3(X_1), u \in \S_1 $ and $\ss \in 
S^{\infty}$. We extend $\ti Z$ to $(\A_{2,loc}(X) 
\times \Ga^+_{2,loc}
(X)) \times S^{\infty}$  by setting
$$
\ti Z(u, \ss) = \ti Z(u|_{X_1}, \ss).
$$
By the construction of $\bold \Psi$, $\ti Z$ has 
no zeros.
\enddefinition

\proclaim{Lemma 3.8} The vector field $\ti Z$ is 
gauge equivariant. 
\endproclaim
\demo{Proof} For $g=g_0 g_1$ with $g_0$ a based gauge
and $g_1 \in S^1$, we have
$$
\ti Z(g^*(u, \ss))= \ti Z(g^*u, g_1^{-1}\ss)=
g_0^{-1} \ti Z(u, g_1^{-1} \ss) $$
$$= g_0^{-1} g_1^{-1} \ti Z(u, \ss)
= g^* \ti Z(u, \ss).
$$
\qed
\enddemo

Now  we replace $Z_e$  by $\ti Z$ in the stable transition
flow equation (3.1). The equation  remains invariant under the 
action of $\G_{2, loc}$ and admits no reducible 
solution. Restricting to 
$S^{2n-1} \subset S^{\infty}$ we 
 obtain the associated 
stablized moduli spaces $\underline
{\MM}_T^0(S_{\a_-} \times S^{2n-1}, S_{\a_+} \times S^{2n-1})$.
Their direct limits as $n \to \infty$ are
 the total moduli spaces 
$\underline \MM_T^0(S_{\a_-} \times S^{\infty}, S_{\a_+}
\times S^{\infty})$. Note that in 
the case $\SS = S^{\infty}$, it suffices for our purpose  
to 
utilize the subspaces $\underline 
\MM_T^0(S_{\a_-} \times S^{2n+1},
S_{\a_+} \times S^{2n+1})$ instead of 
the total spaces $\underline \MM^0_T(S_{\a-} 
\times S^{\infty}, S_{\a+} \times S^{\infty})$.

\proclaim{Theorem 3.9}  The moduli 
spaces $\underline \MM_T^0(S_{\a_-} \times
S^{2n+1}, S_{\a_+} \times S^{2n+1})$ are 
compact Hausdorff spaces. 
Moreover, for generic $ b_0$, we have for 
all $\a_- \in \R_-, \a_+ \in \R_+$  and $n$ 

\noindent (1) $\underline {\MM}_T^0(S_{\a_-} 
\times S^{2n-1}, S_{\a_+}
\times S^{2n-1})$ has the structure of $d$-dimensional smooth
oriented manifolds with corners, where $d=\mu_-(\a_-)
-\mu_+(\a_+) + m_0 - {dim} ~G_{\a_-}+2n$.

\noindent (2) This structure is compatible with the 
natural stratification
of $\underline {\MM}_T^0(S_{\a_-} \times S^{2n-1},$ 
$S_{\a_+} \times S^{2n-1})$.

\noindent (3) The  projections $\pi_{\pm}: 
\underline {\MM}_T^0(S_{\a_-} \times
S^{2n-1}, S_{\a_+} \times S^{2n-1}) \to 
S_{\a_{\pm}} \times S^{2n-1}$ are $S^1$-equivariant
smooth maps.
\endproclaim

Passing to $S^1$ quotient, we obtain the quotient 
version of Theorem 3.9.

With these preparations, it is 
clear  that we can carry over the above 
proof of Theorem 3.1 to the general
case $\SS=S^{2n-1}$  or $\SS=S^{\infty}$.

\head{4. de Rham version}
\endhead

In this section we   construct the de Rham version of equivariant 
Seiberg-Witten Floer homology and cohomology and prove that it is 
isomorphic to the singular version with real coefficients.
We also construct the de Rham 
version of Bott-type Seiberg-Witten Floer homology and cohomology and 
prove that it is isomorphic to the singular version of Bott-type 
theory with real coefficients.
 
\bigskip 
\noindent {\bf Equivariant Theory}
\smallskip

As before, $\SS=S^{2n-1}$ or $\SS=S^{\infty}$. Choose an 
$S^1$ invariant Riemannian 
metric on each $S_{\a}$ and use the 
standard metric on $S^{2n-1}$.  (Of course,
there is no need 
for choosing a metric on the reducible $S_{\a}$.)
Then we obtain the adjoint $d^*$ of the exterior 
differential $d$ on $S_{\a} \times_{S^1} S^{2n-1}$.  
Set $\O^*(\R^0 \times_{S^1} S^{2n-1}) =
\oplus_{\a} \O^*(S_{\a} \times_{S^1} S^{2n-1}).$
We have natural inclusions 
$\O^*(\R^0 \times_{S^1} S^{2n-1} )\subset \O^*(\R^0 
\times_{S^1} S^{2n+1})$ which commute with $d$ 
and $d^*$.   Using this fact we define 
$(\O^*(\R^0 \times_{S^1} \SS), d, d^*)$ to be the direct 
limit of  $(\O^*(\R^0  \times_{S^1} S^{2n-1}), d, d^*)$ as $n$ 
approaches infinity. 

Now we introduce our 
equivariant de Rham complex $C^*_{equ, deR} $,
$$
C^k_{equ, deR} = \oplus_{i+j=k} \O^j(S_i 
\times_{S^1} \SS).
$$
We define $C_*^{equ, deR}$ to be the 
same as $C^*_{equ, deR}$. Next we define the 
desired coboundary operator. We first 
consider the case $\SS=S^{2n-1}$.
For each pair $\a, \b \in \R$ with $\mu(\a) > 
\mu(\b)$ we define a coboundary operator 
$\pa_{\a, \b}^*$ as follows.
If the moduli space $\underline
\MM_T^0(S_{\a} \times_{S^1} S^{2n-1}, S_{\b}
\times_{S^1} S^{2n-1})$ is empty, we define 
$\pa_{\a, \b}^*$ to be the
zero operator. If it is nonempty,  we set
 $\pa_{\a, \b}^* \omega =0$ for
$\omega \not \in \O^*(S_{\b} \times_{S^1} S^{2n-1})$. For $\omega 
\in \O^j(S_{\b} \times_{S^1} S^{2n-1})$,
we set
$$
\pa_{\a, \b}^* \omega =  (-1)^j (\pi_-)_* \pi_+^* \omega,
$$
where $(\pi_-)_*$ denotes the fiber integration 
on 
$\underline {\MM}_T^0(S_{\a} \times_{S^1}
S^{2n-1}, S_{\b} \times_{S^1}
S^{2n-1})$ with respect to the
fibration $\pi_-$. 
(For the simple definition of fiber 
integration see e.g. $\AustinBraamII$.) 
We define the coboundary operator of  the equivariant 
de Rham complex to be
$$
\pa_{equ, deR}^* = d + \sum_{\mu(\a)>\mu(\b)} 
\pa_{\a, \b}^*.
$$

Following  the arguments in $\AustinBraamII$ 
we deduce  $(\pa_{equ, deR}^*)^2=0$.

Since $\O^j(S_i \times_{S^1} S^{\infty})$ is 
the direct limit of $\O^j(S_i \times_{S^1} 
S^{2n-1})$, the above definition of 
$\pa_{equ, deR}^*$ immediately
extends to cover the case $\SS=S^{\infty}$. 
Namely each $\omega \in \O^j(S_{\b} \times_{S^1} 
S^{\infty})$ can be viewed as in  $\O^j(S_{\b} \times_{S^1}
S^{2n-1})$ for some $n$, and hence 
the above definition of the coboundary 
operator can be applied. On the other hand, we 
can  take the  direct limit of the fiber integration
$(\pi_-)_*$ on $\underline \MM_T^0(S_{\a} \times_{S^1}
S^{2n-1},
S_{\b} \times_{S^1} S^{2n-1})$ to obtain a fiber 
integration $(\pi_-)_*$ on $\underline \MM_T^0(S_{\a}
\times_{S^1} S^{\infty}, S_{\b} \times_{S^1} S^{\infty})$.
Then we carry over the above definition of 
$\pa_{equ, deR}^*$ word for word to $\SS=S^{\infty}$.

\proclaim{Definition 4.1} The de Rham version
of equivariant Seiberg-Witten Floer cohomology 
$FH^{SW*}_{equ, deR}$ is defined to be the cohomology of the 
complex $$(C^*_{equ, deR}, \pa_{equ, deR}^*).$$
\endproclaim

Next we construct the  de Rham version
homology theory, using the complex $C_*^{equ, deR}$.
The associated boundary operator is defined as follows.
For $\a, \b \in \R$ with $\mu(\a) > \mu(\b)$ 
and $ \omega \in \O^j( S_{\a} \times_{S^1} \SS)$, we set
$$\pa_{\a, \b} \omega = (-1)^j(\pi_+)_* 
\pi_-^* \omega.
$$
(As before, this operator is defined to be the zero operator for
other $\omega$ or for empty $\underline {\MM}_T^0(S_{\a}
\times_{S^1} \SS,
S_{\b} \times_{S^1} \SS).$)

We define
$$
\pa_{equ, deR} = d^* + \sum_{\mu(\a) > \mu(\b)}
\pa_{\a, \b}.
$$
We have $\pa_{equ, deR}^2 =0.$

\proclaim{Definition 4.2} The de Rham version of
equivariant Seiberg-Witten Floer homology 
$$FH^{SW}_{equ, deR*}$$
is defined to be 
the homology of the complex $ (C^{equ, deR}_*, 
\pa_{equ, deR})$.
\endproclaim

We have the equivariant de Rham index 
filtration  $\F^*_{equ, deR}$ and its
dual filtration $
\F_*^{equ, deR}$,
$$
\F^k_{equ, deR} = \oplus_{j \geq k} 
 \O^*(S_j \times_{S^1}
\SS),
$$
$$\F_k^{equ, deR} = \oplus_{j \leq k}
\O^*(S_j \times_{S^1} \SS).
$$

\proclaim{Theorem 4.3} The equivariant de Rham index 
filtrations  induce spectral sequences 
$E^*_{**}(equ, deR)$
and $E^*_{**}(equ, deR; dual)$ such that the former 
converges to $$FH^{SW*}_{equ, deR},$$ and the latter converges
to $$FH^{SW}_{equ, deR*}. $$ Moreover, we have
$$
E^1_{ij}(equ, deR) \cong H_{deR}^j(S_i \times_{S^1} 
\SS), \tag 4.1
$$
$$
E^1_{ij}(equ, deR; dual) \cong H^{deR}_j(S_i \times_{S^1}
\SS).
\tag 4.2
$$
\endproclaim

This is similar to Theorem 2.2.

\proclaim{Theorem 4.3} The de Rham version of equivariant
Seiberg-Witten Floer homology and cohomology are 
diffeomorphism invariants. Indeed, we have
$(FH^{SW}_{equ , deR})^j \cong (FH^{SW}_{equ; \Bbb R})^j,
(FH^{SW}_{equ, deR})_j \cong (FH^{SW}_{equ;
\Bbb R})_j.$
\endproclaim

\demo{Proof} We present the arguments for cohomologies.
Homologies can be handled in a similar 
way. We have a natural cochain map 
$F: C^*_{equ, deR} \to C^*_{equ; \Bbb R}$, which is 
defined in terms of integrating forms along generalized 
cubical singular 
chains. The fact that this is indeed a cochain 
map follows quite easily from the construction of 
the coboundary operators $\pa_{equ}^*, \pa_{equ, 
deR}^*$ and Stokes' theorem. 

It is easy to see that $F$ preserves the 
index filtrations. The induced homomorphism 
$F^*: GC^*_{equ, deR} \to GC^*_{equ; \Bbb R}$ consists
of the cochain homorphisms from $(\O^*(S_j  \times_{S^1}
\SS), d)$ to $( C^*( S_j \times_{S^1} \SS), \pm \pa_O^*)$
defined in terms of integration. Here $\pm \pa_O^*$
 denotes 
the singular coboundary operator $\pa_O^*$
modified by signs, with 
sign convention $\pm = (-1)^{i+j}$ on 
$ C^i(S_j 
\times_{S^1} \SS)$, cf. [$\WangYe$, Section 7]. 
It follows that the induced cochain homomorphism between 
the $E^1$ terms $E^1_{**}(equ, deR)$ and 
$E^1_{**}(equ, dual;
\Bbb R)$ is an isomorphism. Since we are 
using the real coefficients, this implies, on account of 
the convergence of the involved spectral sequences, the 
desired isomorphism. 
\qed
\enddemo

\bigskip
\noindent {\bf Bott-type Theory}
\smallskip

The construction is parallel to the one above. We have in 
the first place the de Rham version Bott-type 
complex
$$
C^k_{Bott, deR}=\oplus_{i+j=k} \O^j(S_i).
$$
The coboundary operator $\pa^*_{Bott, deR}$
for this cochain complex is defined analogously
to $\pa^*_{equ, deR}$, where we use 
$\underline {\MM}_T^0(S_{\a}, S_{\b})$ instead of 
$\underline {\MM}_T^0(S_{\a}\times_{S^1}
\SS, S_{\b} \times_{S^1} \SS)$. 

\definition{Definition 4.4} We define the de Rham version
Bott-type Seiberg-Witten Floer cohomology $FH^{SW*}_
{Bott, deR}$ to be the cohomology of 
$(C^*_{Bott, deR}, \pa_{Bott, deR}^*)$. We also have 
the de Rham version Bott-type homology 
$FH^{SW}_{Bott, deR*}$.
\enddefinition

We omit the statement of the spectral sequence theorem 
here. In analogy  to Theorem 4.3, we have

\proclaim{Theorem 4.5} We have natural isomorphisms $$FH^{
SW*}_{Bott, deR} \cong FH^{SW*}_{Bott; \Bbb R}, 
FH^{SW}_{Bott, deR*} \cong FH^{SW}_{Bott; \Bbb R*}.$$
\endproclaim

\head{5. Cartan Version}
\endhead

In this section we give a brief account of the 
construction of the Cartan version of equivariant Seiberg-Witten Floer 
homology and cohomology and prove that it is isomorphic to the 
de Rham version.  The construction is in spirit similar 
to that in $\AustinBraamI$. Such a construction 
first appeared in the paper
 $\MarcolliWang$.

We follow the presentation of the Cartan theory given in
$\BGV$. For a compact Lie group $G$  acting smoothly on
a manifold $N$, we have the following space of
$G$-equivariant differential forms on $G$
$$\O_G^*(N)=(\O^*(N)\otimes \Bbb C[{\frak g^*}])^G,$$
which is the subalgebra of $G$-equivariant elements in
the algebra $\O^*(N)\otimes \Bbb C[{\frak g}]$.
Here $\Bbb C[{\frak g^*}]$ denotes the complex polynomial algebra on
the dual Lie algebra $\frak g^*$ of $G$ (i.e. the algebra of
complex valued polynomial functions on $\frak g$), $G$ acts on the ordinary form part by
pullback, and on the Lie algebra part by the adjoint representation.
 The algebra $\O^*(N)\otimes \Bbb C[{\frak g^*}]$
has a $\Bbb Z$-grading defined by
$$\deg(w\otimes z)=\deg(w)+ 2\deg(z),$$
for $w\in \O^*(N)$ and $z\in \Bbb C[{\frak g^*}]$.
There is a natural  differential $d_G$ on $\O_G^*(N)$
$$(d_G\a)(Z)=d(\a(Z))-\iota(Z)(\a(Z)),$$
for $\a\in\O_G^*(N)$, where
$\iota(Z)$ denotes contraction by the 
vector field $Z_N$ induced
by an element $Z\in{\frak g}$ via the action of $G$. One easily verifies
$d_G^2=0$.

Now consider our set-up of the Seiberg-Witten theory on the 3-manifold
$Y$ as before. We introduce the complex $C^*_{Cartan}$,
$$C^k_{Cartan} = \oplus_{i+j=k} \O_{S^1}^j(S_i).$$
For each pair $\a, \b$ with $\mu(\a) > \mu(\b)$, we
define an operator $\pa_{\a, \b}^* : C^k_{Cartan} \to
C^{k+1}_{Cartan}$ as follows. If 
the moduli space $\underline
\MM_T^0(S_{\a}, S_{\b})$ is empty, we 
define $\pa_{\a, \b}^*$ to be the
zero operator. If it is nonempty,  we set 
$\pa_{\a, \b}^* \omega =0$ for
$\omega \not \in \O^*_{S^1}(S_{\b})$. For $\omega \in \O_{S^1}^j(S_{\b})$,
we set
$$
\pa_{\a, \b}^* \omega =  (-1)^j (\pi_-)_* \pi_+^* \omega.
$$
The gauge equivariance of 
the projections $\pi_{\pm}$ of 
$\underline \MM_T^0(S_{\a}, S_{\b})$ ensures that 
this indeed produces an equivariant form. 
We define the coboundary operator of  the Cartan complex to be
$$
\pa_{Cartan}^* = d_{S^1} + 
\sum_{\mu(\a)>\mu(\b)} \pa_{\a, \b}^*.
$$

Following the arguments in $\AustinBraamII$ 
we deduce  $(\pa_{Cartan}^*)^2 =0$.

\proclaim{Definition 5.1} The Cartan version
of equivariant Seiberg-Witten Floer cohomology 
$FH^{SW*}_{Cartan}$ is defined to be the 
cohomology of the complex $(C^*_{Cartan}, 
\pa_{Cartan}^*)$.
\endproclaim

In analogy to the dual de Rham construction, we also have 
a dual Cartan construction which produces  the
homology $FH^{SW}_{Cartan*}$.
We omit the simple details. 

\proclaim{Theorem 5.2} The Cartan version of 
equivariant Seiberg-Witten Floer cohomology 
and homology are diffeomorphism invariants.
Indeed, we have
$FH^{SW*}_{Cartan} \cong FH^{SW*}_{equ,
deR;S^{\infty}}, FH^{SW}_{Cartan*} 
\cong FH^{SW}_{equ, deR;
S^{\infty}*}$.
\endproclaim

\demo{Proof} We follow the arguments in $\AtiyahBott$ and
$\AustinBraamII$.  Consider the principal $S^1$-bundle 
$\pi: S^{2n-1} \to \Bbb C P^n$. We choose the natural 
$S^1$-connection on this bundle which is induced 
by the Euclidean metric structure of $\Bbb R^{2n}$.
Denote its curvature form by $\O(n)$. Using the 
inclusions $S^{2n-1} \subset S^{2n+1}$ we then obtain 
the direct limit of $\O(n)$ as $n$ approaches infinity,
which will be denoted 
by $\O$. Note that $\O$ belongs to 
$\O^2(S^{\infty})$.

Now consider the principal $S^1$-bundles $\pi:
S_{\a} \times 
S^{\infty} \to S_{\a} \times_{S^1} S^{\infty}$ 
with $\Bbb C P^{\infty}$ as fiber. We define 
a natural homomorphism
$$
\ti F: \O_{S^1}^* (S_{\a})\to \O^*_{basic}(S_{\a} \times
S^{\infty}),
$$
where basic forms mean forms which are horizontal and 
$S^1$-invariant.  For generators $\omega \otimes 
u \in \O_{S^1}^*
(S_{\a})$, we set 
$$
\ti F(\omega \otimes u) =  \omega \wedge u(\O).
$$
Composing $\ti F$ with the inverse of the 
natural isomorphism 
$\pi^*: \O^*(S_{\a} \times_{S^1} S^{\infty}) \to
\O^*_{basic}(S_{\a} \times S^{\infty})$ we then 
obtain a natural cochain homomorphism 
$$F_{\a}: (\O_{S^1}^*(S_{\a}), d_{S^1})
 \to (\O^*(S_{\a} \times_{S^1} 
S^{\infty}), d).$$
By $\AtiyahBott$, this homomorphism 
induces an isomorphism $F^*_{\a}$ on cohomologies.

Putting all $F_{\a}$ together we       
obtain $F: C^*_{Cartan} \to C^*_{equ, deR}$
with $F: C^k_{Cartan} \to C^k_{equ, deR}$.
By the constructions of the involved 
coboundary operators it is  easy to see that 
$F$ is a cochain map. 

By the construction, $F$ preserves the 
index filtrations. The induced homomorphism 
$F^*: GC^*_{Cartan} \to GC^*_{equ, deR}$ 
precisely consists
of the cochain homomorphisms  
$F_{\a}: (\O^*_{S^1}(S_{\a}),
d_{S^1}) \to (\O^*( S_{\a} \times_{S^1} S), d)$.
It follows that the induced cochain homomorphism between 
the $E^1$ terms $E^1_{**}(Cartan)$ and $E^1_{**}(equ, deR)$ 
is an isomorphism. On account of 
the convergence of the involved spectral sequences, this
implies the 
desired isomorphism. 
\qed
\enddemo

\head{6. Exponential convergence}
\endhead

In this section we first prove the exponential asymptotics
 result
[$\WangYe$, Propostion 4.2]. We restate it in
a strengthened form  below.

\proclaim{Propostion 6.1}
 Consider a compact set $K_{h, \lambda, \pi}$ of triples 
 $(h, \lambda, \pi)$, where $h$ denotes a metric on 
 $Y$ and $(\lambda, \pi)$ a pair of $Y-$generic parameters 
 corresponding to 
 $h$.    Then there are 
positive constants $\d_0$ and $E_0$ depending only
on $K_{h, \lambda, \pi}$ 
with the following properties.  Let
$(h, \lambda, \pi) \in K_{h, \lambda, \pi}$.
Let  $u=(A,\P)=(a, \p)$ be a temporal trajectory 
 of local $(2, 2)$-Sobolev class  and 
 finite energy.
Then there exist a gauge $g\in \G_2(Y)$ and two smooth solutions
$u_-,u_+$ of (1.1)  such that $\ti u= g^*u$ is smooth
and the following holds. For all $l$,
$$
\|\ti u-u_{+}\|_{l, \d_0; Y \times [T_1, \infty)}  \leq C(l),
$$
$$
\|\ti u - u_-\|_{l, \d_0; Y \times (-\infty, T_2]}
\leq C(l),
$$ 
 where 
$C(l)$ depends only on $K_{h, \lambda, \pi}$, 
$l$ and an upper bound of the energy of $u$, and 
$T_1 >0, T_2 <0$ satisfy $E(u, Y \times [T_1, \infty)),
E(u, Y \times (-\infty, T_2]) \leq E_0$. 
(Recall 
that $E(u,  \Omega) $ means the energy of $u$
 on the domain $\Omega.$)

If $u_-$ and $u_+$ are not gauge equivalent, then
$$
E(u) \geq E_0.
$$

Furthermore, the above weighted Sobolev estimates  
hold for any
smooth 
temporal trajectory $\ti u$ provided that all 
the cited conditions hold and in addition its 
values at $T$  and $T'$ are under smooth control,
with the constants $C(l)$ also depending on 
the said control.

\endproclaim

Note that under the rational homology sphere 
assumption,
the energy $E(u)$ vanishes when $u_-$ and 
$u_+$ are gauge equivalent.  This follows from 
[$\WangYe$, Lemma 2.9].
  
We fix a compact set $K_{h, \lambda, \pi}$. {\it In the
 remainder of 
this section, all the results are under this condition,
and all the constants in the various estimates  depend
on $K_{h, \lambda, 
\pi}$.} 

  Let $u = (A, \Phi)$ satisfy 
the conditions 
in Proposition 6.1. We handle convergence 
at $+\infty$. The situation of $-\infty$ is
similar.

For an arbitary $r$, consider the cylinder $X_{r, r+2} =
Y \times [r, r+2]$. We apply 
[$\WangYe$, Lemma B.1] to produce the Columb form   $u^r$ 
of $u$ over this cylinder. 
By the proof of [$\WangYe$, Lemma 8.6], we have 
$$
\eqalign{
\int_{X_{r, r+2}} |F_A|^2 \leq  &2 \int_{X_{r, r+2}} |F_A 
+ \langle e_i\cdot\p,\p\rangle
*e^i|^2 + 6  \int_{X_{r, r+2}} |\p|^4 \cr 
\leq & 4\int_{X{r, r+2}} |*F_a + \langle e_i \cdot  \p, \p\rangle
*e^i|^2 + 4  \int_{X_{r, r+2}} |\frac{\partial a }{\partial t}
- d_Y f|^2  \cr
+ 6 \int_{X_{r, r+2}} |\p|^4   
\leq & C(1+ E(r)). \cr }
\tag 6.1
$$
   
By [$\WangYe$, Lemma B.1] we then obtain a uniform estimate for 
the $L^{1,2}$ norm of $A-A_0$ over $X_{r, r+2}$, where
 $A_0 = a_0$ 
is a reference connection, with $a_0$ being
a fixed smooth connection over $Y$. 
Now $u^r=(A, \P)=(a+f dt, \p)$ satisfies the 
following elliptic system 
on $Y \times (r, r+2)$  (cf. [$\WangYe$, (2.10)])
$$
\cases &F_A^+ = \frac{1}{4}\langle e_i e_j \P,  \P
\rangle e^i \wedge e^j + \nabla H(a) \wedge dt
+ *(\nabla H(a) \wedge dt), \cr 
&d^*(A-A_0)= 0, \cr
&D_A \P = - \lambda \frac{d}{dt} \cdot \P. \endcases
\tag 6.2
$$

Using this system, the estimate (6.1), 
[$\WangYe$, Lemma 
8.5]  and basic elliptic theory, we then 
derive uniform (independent of $r$) interior smooth 
estimates for $u^r$. 
For an arbitary  sequence $r_k \to \infty$, we obtain
a subsequence, which we still denote by $r_k$, and
an $L^{1, 2} $ weak limit $u_{\infty}=(a_{\infty} + f_{\infty} dt,
 \p_{\infty}
) $ of $u_k := u^{r_k}(\cdot, \cdot -r_k)$, such that 
$f_{\infty} \in L^{1, 2}_0(Y \times [0, 2])$, i.e. $f_{\infty}$ has zero boundary 
trace. (The zero boundary trace condition follows from the Columb
boundary condition of [$\WangYe$, Lemma B.1].) Moreover, 
$u_k$ converges to $u_{\infty}$ smoothly in the interior.  
By the finite energy condition, $u_{\infty}$ satisfies the following
equations (we omit the subscripts $(\lambda, H)$ in 
the perturbed Seiberg-Witten equation):
$$
\hbox{\bf sw}(a_{\infty}(\cdot, t), \p_{\infty}(\cdot, t)) = 0,
\tag  6.3
$$
for all $t \in (0, 2)$, and 
$$
\frac{\partial a_{\infty}}{\partial t} - d_Yf_{\infty}
=0, \frac{\partial f_{\infty}}{\partial t}
-d_Y^*(a_{\infty}-a_0) =0, \frac{\partial \p_{\infty}}{\partial t} +
f_{\infty} \p_{\infty} =0,
$$
on $Y \times (0, 2)$.
It follows that
$f_{\infty}$ is harmonic and hence
$f_{\infty} \equiv 0$. Consequently, 
$u_{\infty}$
has no $t-$component and is $t$-independent. Moreover, 
we have 
$$d_Y^*(
a_{\infty} - a_0 )= 0.
\tag 6.4
$$

We write $u_k = (a_k + f_k dt, \p_k)$ as usual.  
Applying the gauge fixing lemma [$\WangYe$, Lemma 3.2]
 we 
convert $(a_k(\cdot, t), \p_k(\cdot, t))$ for large $k$ and $t$ strickly 
inside $(0, 2)$ into 
satisfying the  gauge fixing condition
$$
d_Y^*( a- a_{\infty}) + \hbox{Im }\langle \p_{\infty}, 
\p \rangle =0.
\tag 6.5
$$

By [$\WangYe$, Lemma 3.2] and basic elliptic estimates we 
can control the involved 
gauges, and hence 
derive that  the new $(a_k, \p_k)$  still
converges to $u_{\infty}$
(we retain the notation $(a_k, \p_k)$ for convenience).

\proclaim{Lemma 6.2} For large $k$ we have 
$$
\|\hbox{\bf sw}( a_k,  \p_k)\|_{L^{2}}^2 \geq c
 \|( b_k, 
\psi_k)\|_{L^{1, 2}}^2
\tag 6.6
$$
for a positive constant $c$, where $ b_k = a_k - a_{\infty},
  \psi_k =
 \p_k - \p_{\infty}.$
\endproclaim
\demo{Proof}  We first express  
$$\sw( a_k,  \p_k) = 
\sw( a_k,  \p_k)- \sw(a_{\infty}, \p_{\infty})$$
 in terms 
of $L_k ( b_k,  \psi_k)$, 
where $L_k $ denotes the linearized Seiberg-Witten operator 
$d\sw(a_{\infty}+ s ( a_k -
a_{\infty}), \p_{\infty} +s( \p_k- \p_{\infty}))$ for
 $s \in [0, 1]$. Then we write  $L_k = 
 d \sw( a_{\infty}, \p_{\infty}) + (L_k-
 d \sw (a_{\infty}, \p_{\infty}))$. By the 
  convergence of $ u_k$ to $u_{\infty}$, 
  the term involving $L_k - d \sw (a_{\infty}, 
  \p_{\infty})$ approaches zero as $k$ goes to 
  infinity.  On the other hand, by the genericity 
  assumption on the perturbation parameters   we have 
 the estimate
 $$
 \eqalign \|d\sw(a_{\infty}, \p_{\infty}) (b', \psi') \|_{L^2}
 \geq c  \|(b', \psi')\|_{L^{1, 2} }
 \tag 6.7
 $$
  for a positive constant $c$ and 
   all $(b', \psi')$ satisfying the gauge condition 
  (6.5).  This is obvious in the case that $u_{\infty}$ 
  is irreducible, for in that case ${ker}~
  d\sw(u_{\infty})|_{G^*_{Y, u_{\infty}}}$ is obviously
  trivial under the genericity assumption. (See 
  [$\WangYe$, Section 3] for the infinitesimal 
  gauge action operator $G_{Y, u_{\infty}}$
  and other relevant materials.)
  If $u_{\infty}$ is reducible, then we have 
  (see [$\WangYe$] for the definition of $\DD_1$)
  $$
  \|\DD_1|_{(a_{\infty}, \p_{\infty}, c) }
  (b, \psi, f)\|_{L^2} 
  \geq  c\|(b, \psi, f)\|_{L^{1,2}},
  $$
  for all $(b, \psi, f) \in {ker}~d_1^* \oplus 
  \Ga_1(S) \oplus \O^0_1(Y)$ which are 
  orthogonal to the kernel of
  $\DD_1$. But this kernel obviously contains $\{0\}
  \oplus \{0\} \oplus i\Bbb R$. Since $\DD_1$ is surjective by 
  the genericity assumption and ${ind}~\DD_1 =1$,
  $\{0\} \oplus \{0\} \oplus \Bbb  R $ must be the 
  entire kernel. It follows that 
  all elements of the form $(b', \psi', 0)$ in
  ${ker}~d^*_1 \oplus \Ga_1(S) \oplus \O^0_1(Y)$
  are orthogonal to ${ker}~\DD_1$. Finally, 
  we observe $\DD_1(b', \psi', 0) =d\sw(a_{\infty},
  \p_{\infty})(b', \psi')$.
  \qed   
\enddemo
              
It is easy to extend the above results to 
a more general setting as follows.
\proclaim{Lemma 6.3} For each 
$\delta >0$  there  is a positive 
number $E_0=E_0(\delta)$ with the following property.
Let $u$ be a solution of 
(1.2) over $Y \times [T_1, T_2]$ with $1 \leq T_2-T_1 
\leq 2$. If $E(u) \leq E_0(\delta)$,
 then 
there is  a solution $ u_{\infty}  = ( a_{\infty}, 
\p_{\infty})$ of (1.1) 
such that $d^*(a_{\infty}-a_0)=0$ and  
$$
\|\ti u-u_{\infty}\|_{C^2(Y \times [T_1+1/2, T_2-1/2]} 
\leq \delta,
$$
$$
\|\ti u-u_{\infty}\|_{C^j(Y \times [T_1+1/2, T_2-1/2]} 
\leq \epsilon(\delta,j),
$$
where $\ti u$ is obtained from $u$ by a suitable 
gauge transformation ($u_{\infty}$ is the same for all $j$),
and the functions $\epsilon(\delta, j)$ approach 
zero as $\delta \to 0$.
If we write $\ti u= (\ti a + \ti f dt, \ti p)$,
then $(\ti a, \ti p)$ satifies the gauge condition (6.5).
Moreover, there holds
$$
\|\sw(\ti u )\|_{L^2}\geq c\|\ti u-u_{\infty}\|_{L^{1, 2}(Y)}
$$
for a positive constant $c$, at each 
$t \in [T_1+1/2, T_2- 1/2].$
\endproclaim

For different sequences $u_k$ as above, we 
may a priori get different limits $u_{\infty}$. 
 However, by the 
monotonicity of 
$\cs$, the limit value 
$\cs(\infty)= \cs(u_{\infty})$ is independent of 
the choice of the sequence. Indeed, it  is equal 
to ${lim}_{t \to \infty}
\cs(a, \p).$  

Now we proceed to derive energy decay. We set $J(t) = 
\cs(t)-\cs(\infty) \equiv \cs((a(t), \p(t))-\cs(\infty).$ 
Working in temporal gauge we deduce 

$$
\frac{\pa J}{\pa t} = \frac{\pa \cs}{\pa t} =
\nabla \cs \cdot \frac{\pa u}{\pa t} 
= -\|\sw (a, \p)\|^2_{L^2}.
\tag 6.8
$$

On the other hand, we have by [$\WangYe$, Lemma 2.5]
$$
J(t) =  \frac{1}{2} E(u, Y \times [t, \infty)).
\tag 6.9
$$
By gauge invariance, the formulas (6.8) and (6.9) 
 hold in any gauge.  
By the above convergence argument, the 
distance from the suitable gauge form of $(a(t), 
\p(t)) := (a(\cdot, t),
\p(\cdot, t)) $ to some 
$u_{\infty}$ goes to zero as $t$ goes to 
infinity.  For each large 
$t$, we 
use the corresponding $u_{\infty}$ and 
the associated gauge. By Lemma 6.2 or Lemma 6.3
we deduce
$$
\|\sw(a(t), \p(t))\|_{L^2} \geq c 
\|(b, \psi)\|_{L^{1,2}}
\tag 6.10
$$
with $(b, \psi) = (a(t)-a_{\infty}, \p(t)- u_{\infty})$.
On the other hand, we have on account of a priori 
estimates in the Columb gauge
$$
\|\sw(s(a(t), \p(t)) + (1-s) (a_{\infty}, 
\p_{\infty}) (b', \psi')\|_{L^2} \leq C \|(b',
\psi')\|_{L^{1,2}}
$$
with $s \in [0, 1]$. Here and in the sequel, $C$ 
denotes a  positive constant whose value depends 
on each context.
Consequently,
$$
\eqalign{ J(t)= & \cs(t)- \cs(\infty) \cr
=&  \int_0^1\sw(s (a(t), \p(t))
+ (1-s) (a_{\infty}, \p_{\infty}))(b(t), \psi(t)) ds
\cr
\leq  & C \|(b(t), \p(t))\|_{L^{1,2}}. \cr }
\tag 6.11
$$
Combining (6.8), (6.10) and (6.11) we then infer

\proclaim{Lemma 6.4} There are positive constants 
$E_0>0$ and $c$ 
such that
$$
J(t) \leq -cJ'(t)
\tag 6.12
$$
whenever $E(u, Y \times [t-1, \infty)) \leq E_0$. 
Consequently,
$$
J(t) \leq J(t_0) e^{-c(t-t_0)}
\tag 6.13
$$
for $t \geq t_0$, provided that $E(u, Y \times
[t_0-1, \infty)) \leq E_0$.
\endproclaim

Next we derive an estimate of $\sw(a, \p)$ in terms 
of energy.

\proclaim{Lemma 6.5} There holds 
$$
\|\sw(a, \p)\|_{L^{\infty}(Y \times 
\{t\}})\leq CE(u, Y \times [t-1, t+1])^{1/2}
\tag 6.14
$$
for $u=(a + f dt, \p)$ in any gauge form.
\endproclaim
\demo{Proof} Fix $t_0$. We use the Columb form 
$u_C$ of $u$ over $Y \times [t_0-1, t_0+1]$.
Write $u_C= (a+fdt, \P)$. 
Using the temporal gauge $g= exp(\int_{t_0-1}^t f)$
we convert $u$ into a temporal form $\ti u=
(\ti a, \ti \p)$, 
for which we have smooth control over 
$Y \times [k-1/2, k+1/2]$. 
We set
$$
v = (\ti b, \ti \psi) =\hbox{\bf sw}(\ti a, \ti \p).
$$
There holds 
$$
\cases \frac{\partial v}{\partial t}
 = &d \hbox{\bf sw}(\ti a ,\ti \p)(v),\cr
G^*_{Y, \ti u} v = &0. \endcases
\tag 6.15
$$
Setting $\ti f =0$, $\ti A= \ti 
b + \ti f dt,$ and $ \ti \Psi(\cdot, t) =\ti \psi$ we 
deduce 
$$
\cases d \SW(\ti u) (\ti A, \ti \Psi) =&0, \cr 
G^*_{X, \ti u} (\ti A, \ti \Psi) =& 0. \endcases
\tag 6.16
$$
This elliptic system has been employed  
in [$\WangYe$, Section 4]. 
Its coefficients are given in terms of $\ti u$.
On the other hand, by gauge invariance  we have  
$$
 \| (\ti A, \ti \Psi) \|_{L^2(Y \times [k-1, k+1])} 
 \leq C E(u, Y \times [k-1, k+1]).
 $$
By basic elliptic theory we then 
obtain higher order estimates for $(\ti A, 
\ti \Psi)$.
In particular, we obtain
$$
\|v\|_{L^{\infty}_{Y \times \{t\}}} \leq C
E(u, Y \times [t-1, t+1])^{1/2}.
$$
\qed
\enddemo

\bigskip

\demo{Proof of Proposition 6.1}

 Using 
the Columb gauges and a patching argument,
we can arrange $u$ to be smooth everywhere. Indeed,
we can start converting $u$  into the Columb gauge
on $Y \times [-2, 2]$. We extend the involved  gauge 
 to the entire $Y 
\times \Bbb R$ by interpolating with the identity.
Next we convert $u$  on 
$Y \times [1, 5]$ and $Y 
\times [-5, -1]$ with a Columb gauge $g$. We interpolate
$g$ such that it becomes the identity over $Y
\times [1, 3/2]$ and $Y \times 
[-1, -3/2]$. Then the resulting $u$ will 
be smooth over $Y \times (-5, 5)$. Moreover, we 
have smooth estimates for it in the interior. 
Arguing like this, we can convert $u$ step by step. 
The final $u$ we obtain is no longer in Columb
gauge everywhere because of interpolation of 
the gauges. But it is smooth and we have uniform
smooth estimates for it.

 We denote this smooth form 
 of $u$ by $u_S$. Now 
we convert $u_S$ into  a temporal form 
by using  the temporal gauge $g=exp(\int_{T_0}^tf_S)$
with $T_0$ satisfying $E(Y \times [T_0, \infty)) \leq
E_0.$  We denote it by $u_P$. It is smooth and
its value at $T_0$ is under smooth control.

By  Lemma 6.4  and 
Lemma 6.5 we deduce 
$$
\|\sw(u_P) \|_{L^{\infty}} \leq C e^{-c(t-T_0)}.
\tag 6.17
$$
By integration, we  first derive  exponential $L^{\infty}$
convergence    for $u_P $
because we have
$$
\frac{\pa u_P}{\pa t}
= \sw(u_P).
$$
We can apply  the  elliptic system   (6.16)  to 
$u_P$ and $\sw(u_P)$.
By basic elliptic estimates we then deduce 
$L^{1, l}$ decay  estimate for $\sw(u_P)$ for any $l$.
Integrating, this implies in turn  
exponential $L^{1, l}$ convergence for $u_P$. 
Differentiating the system (6.16) and iterating 
we then deduce 
higher order decay estimates for $\sw(u_P)$ 
and high order exponential convergence for
$u_P$.     The claimed estimates in exponentially
weighted Sobolev norms follow readily.

It is clear that these arguments apply to 
any temporal form $\ti u$ of $u$ as long as
the value of $\ti u$ at $T_0$ is under smooth 
control.
\qed
\enddemo

The following generalization of 
Lemma 6.4 is very important for us.

\proclaim{Lemma 6.6} There are positive 
numbers $E_0$ and $c$ with the following properties.
Let $u$ be a solution of (1.2) over 
$Y \times [T_1-1, T_2+1]$ such that  
$T_2-T_1 \geq  2$ (this assumption is 
purely for the purpose of 
concise formulation) and $E(u, Y \times [T_1-1, T_2+1])
\leq E_0$. Then at least one of the 
following three cases occurs.

\noindent {\it Case 1}  
$$
E(u, Y \times [t, T_2]) \leq E(u) e^{-c(t-T_1)}
$$
for all $t \in [T_1, T_2].$
\smallskip

\noindent {\it Case 2} 
$$
E(u, Y \times [T_1, t]) \leq E(u)e^{-c(T_2-t)}$$
for all $t \in [T_1, T_2]$. 
\smallskip
 
\noindent {\it Case 3}  There is some $T_0 \in
[T_1, T_2]$ such that
$$
E(u, Y \times [t, T_0]) \leq E(u) e^{-c(t-T_1)}
$$
for all $t \in [T_1, T_0]$ and 
$$
E(u, Y \times [T_0, t]) \leq E(u)e^{-c(T_2-t)}
$$
for all $t \in [T_0, T_2]$.
\smallskip

Note that Case 1 and Case 2 can be considered as special
cases of Case 3.
\endproclaim
\demo{Proof} Applying Lemma 6.3 to  
each interval $[t_1, t_2]$ in $[T_1-1, T_2+1]$ 
of length 1, we obtain a limit $u_{\infty}$.
The limits may depend on the interval, but 
the Chern-Simons value $\cs(u_{\infty})$ doesn't.
Indeed, by the $Y$-genericity assumption,
there are only finitely many elements in 
$\R$, the moduli space of gauge classes of Seiberg-Witten 
points. We choose $E_0$ to be smaller than 
$$
\frac{1}{4}min\{|\cs(\a) -\cs(\b)|: \a, \b \in \R,
\a \not = \b\}.
$$
We also choose $E_0$ so small such that 
the estimates in Lemma 6.3 imply 
$$
|\cs(u(\cdot, t)) - \cs(u_{\infty})| < \frac{1}{4}
min\{|\cs(\a)-\cs(\b)|: \a, \b \in 
\R, \a \not = \b\}.
$$
Employing the equation  [$\WangYe$, (2.18)] relating 
the energy to the Chern-Simons functional, we then 
conclude that $\cs(u_{\infty})$ is independent of 
$u_{\infty}$. We denote it by $\cs(\infty)$.

There are three cases to consider.    We set
$\cs(t)= \cs(a(\cdot, t), \p(\cdot, t))$ as before.

\smallskip

\noindent {\it Case 1} $\cs({\infty}) < \cs(T_2).$

In this case, we consider the function $J(t)= 
 \frac{1}{2} E(u, Y \times [t, T_2]) + \cs(T_2)-
 \cs(\infty) = \cs(t)-\cs(\infty)$ and argue as 
 in the proof of Lemma 6.4. We obtain decay of 
 $J$. Since $\cs(T_2)-\cs(\infty) >0$, this implies 
 decay of the energy and leads to Case 1 in the lemma.
\smallskip
 
 \noindent {\it Case 2} $\cs({\infty}) > \cs(T_1).$
 
 This is similar to Case 1, with the roles of $T_1$ and 
 $T_2$ reversed. It corresponds to the  
 negative time infinity version of Lemma 6.4 and 
 leads to Case 2 of the lemma.

\smallskip
 \noindent {\it Case 3} $\cs(T_1) \geq \cs({\infty}) 
 \geq \cs(T_2).$
 
 In this case, there is a time $T_0 \in
 [T_1, T_2]$ such that $\cs(\infty)= \cs(T_0)$.
 We consider $J(t)= \cs(t)-\cs(T_0)$ for 
 $t \in [T_1, T_0]$ and $J(t) = \cs(T_0) 
 -\cs(t)$ for $t \in [T_0, 
 T_2]$, and arrive at  Case 3 of the lemma.
 \qed
 \enddemo 
 
A simple consequence of Lemma 6.6 is the 
following result.

\proclaim{Lemma 6.7} Assume the same as in 
Lemma 6.6. 
Consider the local energy of 
$u$
$$
E_u(t) \equiv  E(u, [t-1/2, t+1/2]) .
$$
At least one of the following three
cases occurs.
\smallskip

\noindent {\it Case 1} 
$$
E_u(t) \leq E(u) e^{-c(t-T_1-1/2)}$$
for all $t \in [T_1+1, T_2-1]$.

\smallskip

\noindent {\it Case 2}

$$
E_u(t) \leq E(u) e^{-c(T_2-t-1/2)}
$$
for all $t \in [T_1+1, T_2-1]$. 
\smallskip

\noindent  {\it Case 3} 
$$
E_u(t) \leq 2E(u) e^{-c(t-T_1-1/2)}
$$
for all $t \in [T_1+1, T_0]$ and 
$$
E_u(t) \leq 2E(u) e^{-c(T_2-t-1/2)}
$$
for all $t \in [T_0, T_2-1]$, where 
$$
T_0 =(T_1+T_2)/2.
$$
\endproclaim

Next we derive a uniform pointwise estimate. 

\proclaim{Lemma  6.8} Let $u$ be a 
temporal trajectory satisfying an energy bound 
such that its value at some time $t_0$ is under 
smooth control. Then we have uniform pointwise 
smooth estimates for $u$. If we only assume 
an energy bound for $u$, we can find 
a gauge $g \in \G_3(Y)$ such that the value of
$g^* u$ at $t=0$ is under smooth control in dependence 
on the energy bound. Consequently, we have smooth estimates 
for $g^*u$. 
\endproclaim 
\demo{Proof} The last statement follows from the 
arguments in the proof of Proposition 6.1.
 Assume 
$E(u) \leq C$ for some $C$ and smooth control over the value of 
$u$ at $t_0$. For convenience, assume $t_0=0$.
(Our arguments work for any $t_0$.)
We write $\Bbb R$ as a union of $m$
intervals $I_j=[t_j, t_{j+1}], t_j < t_{j+1},$
with $m \leq 2E(u)/E_0$, such that
for each $j$, either both $t_{j+1}-t_j > 1$ and 
$$
E(u, [t_j-1, t_{j+1}+1]) \leq E_0
$$
hold,  or 
$t_{j+1}-t_j =1$. Here, we adopt the convention
that e.g. $[t_j, t_{j+1}]$ means $(-\infty, t_{j+1}]$
if $t_j= -\infty$. To derive the desired 
estimates, we start with the interval $I_{j_0}$ containing 
$0$.  If $E(u, I_{j_0}) \leq 
E_0$, we apply Lemma 6.7,  the assumption 
on $u(\cdot, 0)$ and the arguments in the proof of 
Proposition 6.1 to derive uniform pointwise smooth estimates
for $u$ over $Y \times I_{j_0}$. If $t_{j_0+1}-
t_{j_0} =1$, we apply Lemma 6.5 and the arguments 
in the proof of Proposition 6.1 to achieve the same.
 With this, we move on to 
the adjacent intervals of $I_{j_0}$. By induction, we 
obtain the desired estimates for $u$.
\qed 
\enddemo

The next lemma  concerns with an energy estimate.

\proclaim{Lemma 6.9} 
For each positive number $\varepsilon$, there 
exists a positive number $\varepsilon'$  with the following properties.
Consider a  trajectory 
$u$ and an interval $[T_1, T_2]$. Assume 
$$
max\{E_u(t): t \in [T_1, T_2]\} \leq \varepsilon'.
$$
Then
$$
E(u, [T_1, T_2]) \leq \varepsilon.
$$
\endproclaim
\demo{Proof} 
Consider a sequence $u_k$ of trajectories along 
with a sequence of intervals $[T_{k, 1}, T_{k, 2}]$ 
such that 
$$
max\{E_{u_k}(t): t \in [T_{k, 1}, T_{k, 2}]\} \to 0.
\tag 6.18
$$
We claim that $E(u_k, [T_{k, 1}, T_{k, 2}])
\to 0$, which leads to the desired 
conclusion of the lemma. For simplicity, let's assume that  
 $[T_{k, 1}, T_{k, 2}] \subset 
[T_{k+1, 1}, T_{k+1, 2}]$.  The general case can 
be handled by the same argument with slight modification.
 The assumption (6.18) implies that 
in local Columb gauges, the solutions converge 
to Seiberg-Witten points. If these limits are all on 
the same Chern-Simons level, 
then the global energy must converge 
to zero in view of the equation  [$\WangYe$, (2.18)]
relating the
energy to the Chern-Simons functional. Hence, if 
the claim is false,   
there must be at least two different 
Chern-Simons levels for the limits.   By (6.18),
 they must  be obtained at 
two support locations whose time 
distance approaches infinity. Now we consider the intermediate 
limits. At least one of them must lie on 
a Chern-Simons level different from the above 
two levels.  Otherwise we can   easily derive  a contradiction to  
the fact that the above two levels are distinct.

But the time distance from 
   the  support location of this new limit 
   to the support location of the previous two
   limits approaches infinity.   
   Hence  we have rooms in between for finding  
further levels. By the same arguments,
these levels are pairwise different. This leads 
to a contradiction to the generic parameter 
assumption, which implies that there are only 
finitely many Chern-Simons levels.
\qed
\enddemo

\proclaim{Corollary 6.10} There is a positive 
number $E_0^*$ with the following 
properties. Assume that $u \in \N(\a,
\b)$ with $\a \not = \b$. Then 
$$
max\{E_u(t): t \in \Bbb R\} \geq  E_0^*.
$$
\endproclaim

Now we proceed to analyse convergence of
temporal
trajectories of uniformly bounded energy. 
The appropriate concept of convergence is 
smooth  
piecewise exponential convergence 
introduced in 
[$\WangYe$, Definition 6.11].  
They are measured in terms of  the distance 
$d_{l, \bold r}$ given in [$\WangYe$, Definition 6.10]. 
Note that a pair of positive exponents $\d=(\d_-, \d_+)$ are 
involved in this convergence concept.  We shall refer 
to $\d$ as the {\it convergence exponent}.

\proclaim{Theorem  6.11}  Let $u_k$ be a sequence of 
smooth temporal trajectories with uniformly 
bounded energy. We also assume that their energies
are uniformly positive. Moreover, 
  assume that $u_k(\cdot, 0)$ are under  uniform
smooth control. Then there is a subsequence of $u_k$,
still denoted $u_k$, along with a sequence of 
numbers $t_k$ such that $\tau_{-t_k}(u_k) 
\equiv u_k(\cdot, \cdot +t_k)$   converge 
in smooth piecewise exponential fashion to a proper temporal  
piecewise trajectory $u=(u^1,...,u^m)$ with 
$m$ bounded from above by 
a constant depending only on the energy bound.
Here, the convergence exponent depends only on
$K_{h, \pi, \lambda}$.
If $u_k(\cdot, 0)$
are not under uniform smooth control, they can be 
converted into so by a sequence of gauges $g_k \in
\G_3(Y)$.

If $E(u_k)$ are not uniformly positive, then
we have the same conclusions with $t_k=0$ and 
the limit  $u$ being 
a time-independent smooth trajectory (in particular
$m=1$).
\endproclaim
 
\demo{Proof} \hskip 1 cm
Consider a sequence of smooth 
temporal solutions  $u_k$ with uniformly bounded 
energy whose values at $t=0$ are 
under uniform smooth control. Moreover, assume 
that their energies are uniformly positive.
 By Lemma 6.8, we have 
uniform pointwise smooth estimates for 
$u_k$.   Choose $t_k$ such that 
$E_{u_k}(t_k) \geq E_0^*$. We  
deal with  the 
sequence $\tau_{-t_k}(u_k)$ instead of $u_k$.
For convenience, we still denote it by $u_k$.
By convergence
 on compact domains we obtain 
a limit  $u_{\infty,1} $ from $u_k$.  (We pass to 
subsequences whenever necessary.)  We first analyse the 
case that 
 the limit of energy  is entirely captured by 
$u_{\infty, 1}$,  i.e.
$$
E(u_{\infty, 1}) = lim~ E(u_k).
$$
(We may assume that the limit of energy exists.)  Then there 
is a number $T>0$ such that for large $k$,
the energy of $u_k$ outside of the domain $Y \times [-T, T]$ 
is smaller than the critical level $E_0$. By Propositon 6.1 
we obtain uniform exponential estimates for $u_k$,
which, along with the convergence on compact 
domains, imply uniform convergence in exponential 
norms with a smaller exponent. More precisely,
we have
$$
d_{l, \d}(u_k, u_{\infty, 1}) \to 0
$$
for $ \d =(\d_-, \d_+)$ with 
$\d_, \d_+ \in (0, \d_0) $ and all $l$.

This is the ``one piece
" case of piecewise exponential convergence. 

We remark in passing 
that the case of the energies $E(u_k)$ being 
not uniformly positive can be handled in 
a similar way.

Next assume that some energy is lost near infinity, 
i.e.
$$
E(u_{\infty, 1}) < lim~E(u_k).
$$
(Compare $\YeGromov$ where an energy loss 
analysis for pseudo-holomorphic curves is 
presented.) 
Then we can find a sequence of positive 
numbers $T_k $ approaching infinity such that 
the energy of $u_k$ over the domain $Y \times 
[T_k, \infty)$ or the domain $Y \times [-\infty,
-T_k)$ does not converge to zero. 
We handle the first, while the second is similar.

We set $t_k^1=0$. By the energy loss assumption
and Lemma 6.9,
we can find times $t_k^2 \to \infty$ 
such that the local energies $E_{u_k}(t_k^2)$ 
are uniformly positive.
Now we use $-t_k^2$ as the translation amount and 
consider $\ti u_k =\tau_{-t_k^2}( u_k)$.
In other words, we translate the times $t_k^2$ 
back to zero. The 
limit $u_{\infty, 2}$ we get from $\ti u_k$ has 
positive energy, and hence is nontrivial. Indeed, 
$E(u_{\infty, 2})$
 will be at least a quantum portion, namely
it is no smaller than the number $E_0$ given 
in Proposition 6.1.
Now we ask whether some energy might be 
lost between the first limit and the second.
More precisely, we ask whether there are 
times $T_k, T_k'$ with 
$T_k \to \infty, T_k < T_k'$ and 
$t_k^2-T_k' \to \infty$, such that 
$E(u, Y\times [T_k, T_k'])$ 
does not converge to zero.
In view of Lemma 6.9, this is equivalent 
to the question whether there are times 
$t_k^3$ such that $t_k^3 \to \infty, t_k^2-t_k^3 
\to \infty$ and the local energies $E_{u_k}(t_k^3)$
are uniformly positive. 
Assume so. Applying translations in the amount 
$-t_k^3$  
we 
find a nontrivial limit. We say that 
the ``support location" of this limit 
lies between those of the first and 
second limits.  Taking into 
account the time order, we denote this limit by
$u_{\infty, 2}$, while renaming the second limit 
$u_{\infty, 3}$.  

We also ask whether 
there might be lost energy in support locations 
above that of the second limit. We consider 
here instead $t_k^3 \to \infty$ with 
$t_k^3 -t_k^2 \to \infty$. 

Arguing this way, and also taking into account 
 the $-\infty$ direction,
we 
obtain  limits $u_{\infty}^1, u_{\infty}^2, ...$
with increasing time order for support locations. 
Since the energy of each limit is at least $E_0$,
and we obviously have 
$$
\sum E(u_{\infty}^j) \leq lim~ E(u_k),
$$
we can obtain only finitely many limits and arrive at
$$
\sum_{1 \leq j \leq m} E(u_{\infty}^j) = lim~E(u_k)
\tag 6.19
$$
for some $m$.
In other words, all energy loss is accounted for.
By the quantum energy property of the limits, we have 
an upper bound for $m$ in terms of an energy bound.

By the energy identity (6.19), we deduce
$$lim sup_{k \to 
\infty} E(u_k, Y \times [t_k^j +T, t_k^{j+1} -T'])
\to 0,
\tag 6.20
$$
as $T$ and $T'$ both approach infinity. 

Now we analyse the convergence of $u_k$ to the
$m$-tuple $(u_{\infty}^1,...,u_{\infty}^m)$ in more 
details. For simplicity, we consider the case 
$m=2$. The arguments extend to the 
general case straightforwardly.  By the energy 
control (6.20), we obtain  for 
large $T, T'$ Seiberg-Witten points $u_{\infty}$ 
corresponding to  
$\ti u_k$ along $ [T, t_k^2-T']$ as given by 
Lemma 6.3. By the arguments in the proof 
of Lemma 6.6, the Chern-Simons value 
$\cs(\infty)$ is the same for all these 
points. On the other hand, it is easy 
to see (by a simple limit argument) that 
$$
\cs(\infty) = \cs(\ti u_k(\cdot, T_k))
$$
for some $T_k \in (T, t_k^2 -T')$, as long 
as $k$ is big enough.  By Lemma 6.8
we  obtain decay estimates 
$$
E_{\ti u_k}(t) \leq Ce^{-ct}
$$
for $t \in  [0, t_k^2/2]$,
and 
$$
E_{\ti u_k}(t) \leq  Ce^{-c(t_k^2-t)}
$$
for 
$t \in [t_k^2/2, t_k^2],
$
where $C$ depends on an energy bound.
By the arguments in the proof of 
Proposition 6.1, we then deduce 
the exponential convergence of 
$\ti u_k$ on $Y \times (-\infty, t_k^2/2]$ 
to $u^1_{\infty}$ with 
$t_k^2$ seen as in the positive 
direction, and on $Y \times [t_k^2/2, 
+\infty)$ to $u^2_{\infty}$, with $t_k^2/2$ seen as 
in the negative direction.
We also derive for $t, t' >0$  
$$
\|\ti u_k(\cdot, t)-\ti u_k(\cdot, t_k^2-t')\|_{L^{\infty}}
\leq C (e^{-ct}+ e^{-ct'}).
$$
Taking limit, we deduce that 
$$
\|u^1_{\infty}(\cdot, t)-
u^2_{\infty}(\cdot, -t')\|_{L^{\infty}}
\leq C(e^{-ct}+ e^{-ct'}).
$$
This implies that $lim_{t \to +\infty} u^1_{\infty}
= lim_{t \to -\infty} u^2_{\infty}$. Consequently,
$(u_{\infty}^1, u_{\infty}^2)$ is a  
piecewise 
trajectory, cf. [$\WangYe$, Definitions 6.5]. It is proper 
because both $u_{\infty}^1$ and $u_{\infty}^2$ have 
nonzero energy.
 By the  obtained estimates one readily deduces 
 that $u_k$ converges 
to $(u_{\infty}^1, u^2_{\infty})$ in smooth 
piecewise exponential fashion. 
\qed
\enddemo

\demo{Proof of [$\WangYe$, Thereom 6.15]}
The compactness part is an immediate consequence 
of Theorem 6.11.  To show the Hausdorff property of 
$\underline 
\MM_T^0(S_{\a}, S_{\b})$, consider 
$\omega_1=(\omega_1^1,...,
\omega_1^k) \in 
\underline \M_T^0(S_{\a}, S_{\b})_k, \omega_2=
(\omega^1_2, ..., \omega_2^j) 
\in \underline \M_T^0(S_{\a}, S_{\b})_j, \omega_1 
\not  = \omega_2$. Choose 
piecewise trajectories $u_0, v_0$ to  represent 
$\omega_1$ and $\omega_2$ respectively. 
For $r >0, \varepsilon>0$ we set $\bold r_k = (r,..., r) 
\in \Bbb R^k_+$ and define
$$
U_r(\omega_1)= \{ [u]^T_0 \in
\M_T^0(S_{\a}, S_{\b}): \hbox{ There 
is a } u 
\in [u]^T_0 \hbox{ such that } $$
$$
d_{2, \bold r_k}(u, u_0 \# \bold r_k) 
< \varepsilon \},
$$
where $\#$ is the suspension (pre-gluing) map
introduced in [$\WangYe$, Section 6].
Similarly, we have $U_r(\omega_2)$.
For $r$ large enough and $\varepsilon$ 
small enough, the neigborhoods $\bar U_r(\omega_1)$
and $\bar U_r(\omega_2)$ are disjoint.
\qed
\enddemo

The proof of [$\WangYe$, Proposition 6.18] is similar,
we leave the details to the reader.

\head{7. Structures near infinity} 
\endhead

The main purpose of this section is 
to establish the smooth structure of 
the moduli spaces $\underline \MM_T^0(S_{\a},
S_{\b})$ (the temporal model), which is equivalent to the 
smooth structure of the moduli 
spaces $\underline \MM^0(p, q; \S\W_0)$
(the fixed-end model)
with $p \in \a, q \in \b$ (see [$\WangYe$, 
Section 6]). 
We show namely that they are  
smooth manifolds 
with corners. In other words, they are 
modelled on $\bar \Bbb R_+^d$, where $d$ is 
the relevant dimension. (Recall $\Bbb R^+ = (0, \infty)$. )
  There  is 
a natural structure of stratification for these 
moduli spaces, which is provided by 
its subspaces of $k$-trajectories. More precisely,
we have 
$$
\underline \MM^0_T(S_{\a},
S_{\b}) = \cup_k \underline \M^0_T(S_{\a},
S_{\b})_k,
$$
$$
\underline \MM^0(p, q; \S\W_0) =\cup_k \underline 
\M^0(p, q; \SS\W_0)_k,
$$
cf.  
[$\WangYe$, Section 6]. 
To establish the said smooth structure, it remains 
to construct coordinate charts along these strata.
We present a detailed  treatment of the 
temporal model. We also sketch a treament of 
the fixed-end model (independent of the 
treatment of the temporal model). 
Of course,
we only need either one for the constructions of 
our Seiberg-Witten Floer theories.  However, it is 
useful for conceptual understanding to 
clarify the both models.

In the following, we assume that 
the relevant perturbation parameters are 
generi-
c.  
 
\proclaim{Proposition 7.1} Consider (distinct)
$\a, \b \in \R$. 
There are neighborhoods $U_k$, $\hat U_k$ of $\underline 
{\M}^0_T(S_{\a}, S_{\b})_k$
in $\underline {\M}^0_T(S_{\a}, S_{\b})_k\times
 [0, \infty)^{k-1}$ and
$\underline {\MM}^0_T(S_{\a}, S_{\b})$ respectively, and a
homeomorphism $\hbox{\bf F}_k:U_k \to \hat U_k$ such that
the restriction of $\hbox{ \bf F}_k $ 
to $U_k^0=U_k \cap (\underline
{\M}^0(S_{\a}, S_{\b})_k \times (0, \infty)^{k-1})$
is a diffeomorphism. Moreover, the following hold:

(1)  The restriction of 
$\bold F_k$ to $\underline \M_T^0(S_{\a}, 
S_{\b})_k$, which 
is identified with $\M_T^0(S_{\a}, S_{\b})_k$ 
$\times \{(0,...,0)\}$,
is the identity map.
 
(2) For each compact set $K$ in $\underline
 {\M}^0_T(S_{\a}, S_{\b})_k$, there is
a positive number $r_0$ such that
 $K \times [0, r_0]^{k-1} \subset U_k$.

(3)
For $1 \leq j \leq k-1$,
the restriction of $\hbox{\bf F}_k$ to the $j$-th 
open boundary
stratum of $U_k$ is a diffeomorphism onto the $j$-th 
open boundary
stratum $\hat U_k \cap\underline {\M}^0_T(S_{\a}, S_{\b})_{j}$
of $\hat U_k  \cap \underline
{\MM}^0_T(S_{\a}, S_{\b})$. Here e.g. the first
boundary stratum of $\underline
{\M}^0_T(S_{\a}, S_{\b})_k \times
[0, \infty)^{k-1}$ is
$\underline {\M}^0_T(S_{\a}, S_{\b})_k  \times (\{0\} 
\times  (0, \infty)^{k-2}
\cup (0, \infty) \times \{0\} \times (0, \infty)^{k-2}
...\cup (0, \infty)^{k-2} \times \{0\}).$

(4) For different $k, j$,
the transitions between $\bold F_k$ and $\bold F_j$ are 
smooth.

The maps $\hbox{\bf F}_k$ define the structure of
smooth manifolds with corners for 
the temporal model $\underline 
{\MM}^0_T(S_{\a}, S_{\b})$ as
stated in [$\WangYe$, Theorem 6.19].
\endproclaim

Obviously, the crucial result [$\WangYe$, Theorem 6.19]
on the smoothness of $\underline \MM_T^0(S_{\a}, S_{\b})$
follows from this proposition.

The maps $\bold F_k$  will be constructed by a gluing 
process. First, we need to construct suitable 
presentation models for our moduli spaces.

Consider a temporal Seiberg-Witten
trajectory $u 
\in \N_T(S_{\a}, S_{\b})$ (see
[$\WangYe$, Section 6] ) for distinct $\a$ and $\b$.
 We define $\rho_+(u)$ and $\rho_-(u)$
by the following equations
$$
E(u, Y \times [\rho_+(u), \infty))= E_0/2, 
E(u, Y \times (-\infty, \rho_-(u))) = E_0/2,
$$
where $E_0$ is given in Proposition 6.1.
Let $<u>$ denote the set of trajectories which
are gotten from $u$ by a time translation.
Let $u^*_+(u) $ denote the element in $<u>$ whose $\rho_+$
value equals zero, and $u^*_-(u)$ denote the one whose 
$\rho_-$ value equals zero.

\proclaim{Lemma 7.2} The elements $u^*_+(u)$
and $u^*_-(u)$ are uniquely determined. Moreover,
they each give rise to a  
transversal slice for the time translation
action on $\N_T(S_{\a}, S_{\b})$.
\endproclaim
\demo{Proof} First, we show that e.g. 
$u^*_-(u)$ is well-defined. 
It is easy to see that 
there is some $\ti u= \tau_t(u)$ in $<u>$ such that 
$E( \ti u, Y \times [0, \infty)) = E_0/2$.
If there is another, then it follows that 
the energy of $u$ over  
$Y \times I$ for a nontrivial interval $I$ is zero. By 
unique continuation (see 
below for more details),
 $\partial u / \partial t$ 
must be identically zero.
This is impossible, because $\a \not =
\b$.
 
Next we
consider the function $
E(u, [0, \infty))$ on $\N_T(S_{\a}, S_{\b})$.
 We have
$$
 \frac{d}{dt} E( \tau_t(u), [0, \infty))|_{t=0}=
 \frac{d}{dt} E(u, [t, \infty))|_{t=0} = \int_{Y \times 
\{0\}} |\frac{\pa u}{\pa t}|^2 .
$$
If this is zero, then $\pa u / \pa t $ is identically 
zero at $t=0$. Differentiating the equation
$\pa u/ \pa t = \sw(u)$  we derive 
that all time derivatives of $u$ at $t=0$ 
vanish. By unique continuation,  
this implies that ${\partial u}/{ 
\partial t}\equiv 0$, which is impossible.
More precisely,
we consider the quantity $v= \sw(u)$. First, 
we know that $v=0$ at time $t=0$. But 
${\partial v}/{\partial t} = d \sw (v) =0$
at $t=0$. Arguing in this fashion we deduce that 
all space and time dertivatives of $v$ at $t=0$ vanish. 
On the 
other hand, $v$ satisfies the elliptic system
(6.15) or rather (6.16), which 
implies a second order system 
with  scalar symbol. Hence the unique 
continuation principle holds. It follows that $v 
\equiv 0$. 

We conclude that the derivative $d E(\tau_t(u),
[0, \infty)) /dt
\not = 0$. The desired transversality follows.
\qed
\enddemo

Let $\N_T^s(S_{\a}, S_{\b}) $ denote the transversal 
slice given by $u^*_-$, i.e.
$\N_T^s(S_{\a}, S_{\b})$ is the submanifold of 
$\N_T(S_{\a}, S_{\b})$ defined by the equation $ \rho_- =0$.
We have corresponding slices $\N_T^s(S_{\a}, S_{\b})_k, 
\NN_T^s(S_{\a}, S_{\b})$ for the corresponding 
spaces of piecewise trajectories.
By the gauge invariance of energy, 
all these slices are preserved under actions
by gauges in $\G_3(Y)$.  (Because of the temporal
condition, only such gauges are admitted.)

Next we construct a global slice for the 
action of based gauges. Consider the 
global slice $\S_1$ for the action of 
$\G_3^0(Y \times [-1, 1])$ on 
$\A_2(Y \times [-1, 1]) \times 
\Ga^+(Y \times [-1, 1])$ introduced in 
[$\WangYe$, Lemma 8.2]. Using the temporal 
transformation $g_T$ introduced in 
[$\WangYe$, Lemma 6.1] we convert $\S_1$ into
a global slice $\ti \S_1$ consisting of temporal elements. 
We set
$$
\S_T(\a, \b)= \{u \in \N^s_T(S_{\a}, S_{\b}):
u|_{Y \times [-1, 1]} \in \ti \S_1\}.
\tag 7.1
$$ 

The following lemma is readily proved.

\proclaim{Lemma 7.3} The space  
$\S_T(\a, \b)$ is a global transversal slice of the 
action of $\G^0_3(Y)$ on $\N_T^s(S_{\a},
S_{\b})$.  Consequently, it is 
diffeomorphic to the moduli space 
$\underline \M^0_T(S_{\a},
S_{\b})$ and hence can be used as its 
presentation model.
By the construction and 
the arguments at the beginning of Section 6,
 we have smooth control over 
$u(\cdot, 0)$ for $u \in \S_T(\a, \b)$. 
\endproclaim

Next we construct presentation models
for the moduli  spaces 
of consistent multiple temporal 
trajectory classes $\M^0_T(
S_{\a_0}, ..., S_{\a_k})$.
First, we have the fibered 
product 
$$\eqalign{ \S_T({\a_0}, \a_1) 
\times_{S_{\a_1} } 
\S_T(\a_1, \a_2) \times_{S_{\a_2}} 
... \S_T(\a_{k-1}, \a_k)
= \cr 
\{(u_1, ..., u_k) \in
\S_T(\a_0, \a_1) 
\times ... \S_T(\a_{k-1}, \a_k): \cr
\pi_+(u_i) = \pi_-(u_{i+1}),
i=1, ..., k-1\} \cr}
$$
as a presentation 
model.  However, the elements here 
may not be piecewise trajectories, i.e. 
the endpoints of 
$u_i$ as above may not match each other.
We modify this model as follows.
Let $End_{\pm}$ denote 
the endpoint (i.e. limit) at $\pm \infty$.
For $u=(u_1, ..., u_k)$ in this model,
we set $\ti u = (u_1, u_2 
+ End_+(u_1)-End_-(u_2), 
..., u_k + \sum_{1\leq i \leq k-1 } (End_+(u_i)- End_-(u_{i+1}))
$. Obviously, $u$ is a temporal piecewise
trajectory.  By the consistent condition,  the adjustments  
involved in each portion are in terms of actions of 
elements in $\G^0_3(Y)$. Hence we obtain our 
desired  better representation model
for $\underline \M^0_T(S_{\a_0}, ..., 
S_{\a_k})$:
$$
\S_T(\a_0, ..., \a_k) = \{ \ti u:  u \in
\S_T({\a_0}, \a_1) 
\times_{S_{\a_1} } 
\S_T(\a_1, \a_2) \times_{S_{\a_2}} 
... \S_T(\a_{k-1}, \a_k) \}.
$$
 By taking unions we then obtain 
 representation model $\S_T(\a, \b)_k$ for 
 $\underline \M^0_T(S_{\a}, S_{\b})_k$, 
 and  representation model 
 $\SS_T(\a, \b)$ for 
 $\underline \MM^0_T(S_{\a}, S_{\b})$.
 We set 
 $$
 \ti \SS_T(\a, \b) = \cup_{k}(\S_T(\a, \b)_k 
 \times \Bbb R_+^{k-1}).
 $$

The following lemma follows from straightforward 
computations.  
Here, we use again the suspension or 
pre-gluing operator $\#$ introduced in
[$\WangYe$, Section 6]. The piecewise 
exponential weight $w_{u \sharp  \bold r}
= w_{\bold r}$ with $\bold r = 
(r_1, ..., r_{k-1}) \in \Bbb R_+^{k-1}$ 
and the corresponding weighted norm $\|\cdot 
\|_{l, w_{u \# \bold r}}$
were  also introduced in [$\WangYe$, Section 6].
These quantities involve a pair of  exponents $
\d=(\d_-, \d_+)$.
Throughout the sequel we fix a pair $\d$ with 
$\d_-, \d_+ \in (0, \d_0)$,
where $\d_0$ is given in Proposition 6.1.

\proclaim{Lemma 7.4} Let $K$ be a compact
set in $\S_T(\a, \b)_k$.
 Then there are 
constants $C=C(K)$ and $r_0= r_0(K)$
such that  
$$
\|\SW( u \sharp \bold r)\|_{1, w_{u \sharp 
\bold r}} 
\leq C e^{-(\d_0- \d)r_{min}}
\tag 7.2
$$
for $ u\in K, r_i \geq r_0$, where $r_{min} =\min~ \{r_1,...
r_{k-1}\}.
$
\endproclaim

\proclaim{Proposition  7.5} Let $K$ be 
a compact set in $\S_T(\a, \b)_k$.
There are  positive
numbers $\bar r= \bar r(K)$ and
$\sigma = \sigma(K)$ with the following
properties. 
For $u \in K$, $\bold r = (r_1, ...,
r_{k-1} )\in \Bbb R^{k-1}_+$ and $\bar u =(u, \bold r)$,
consider the linearization of $\SW$ at $\sharp
\bar u  \equiv u \sharp \bold r$, $d\SW_{\sharp \bar u}:
%\Omega^0_{1, \d}
\Omega^1_{2, \d} \oplus \Ga^+_{2, \d}\to
\Omega^+_{1, \d}
\oplus \Ga^-_{1, \d} $.
If $$r_i \geq 
\bar r, i=1, ..., k-1, \tag 7.3
$$ 
then
$d\SW_{\sharp \bar u}$ has a right inverse $Q=Q_{\bar u}$
with
$$
\|Q\|_{w_{\bar u}} \leq C= C(K), \tag 7.4
$$
where the subscript
$w_{\bar u}$ means that the norm
is measured in terms of the
$w_{\bar u}$-weighted Sobolev norms
instead of the $\d$-weighted norms, 
cf [$\WangYe$, Section 6].
Furthermore, the equation $\SW(\sharp \bar u +
Q \zeta) =0$ with $\zeta \in \Omega^+_{1, \d}
\oplus \Ga^-_{1, \d}, \|\zeta \|_{1, w_{\bar u}}
\leq \sigma
$ has a unique solution
$\zeta(\bar u)$. 

For $\bar u= (u, \bold r)$ with $r_1, ..., r_{k-1} \geq 
\bar r$ we define
$$
 {\G l}_{K} (\bar u) = \sharp  \bar u + 
Q \zeta(\bar u)
$$
and
$$
\G l^T_K (\bar u)= T_G(\G l_K(\bar u)),
$$
where $T_G$ is the temporal transformation, 
see [$\WangYe$, Section 6].

For $\bar u$ with $r_1, ..., r_{k-1} \in 
(0, 
\bar r^{-1}]$
we define
$$
\ti \G l^T_{K} ( \bar u) =  \G l^T_K (Inv(\bar u)),
$$
where
$$
Inv(\bar u) = (u, (r_1^{-1}, ..., r^{-1}_{k-1})).
$$

These gluing maps for different $K$ coincide with each 
other on 
the intersection of their definition domains,
hence we shall replace the subscript $K$ 
by $k$. Thus we have $\G l_k, \G l_k^T$ and $
\ti \G l_k^T $.
Taking an exhaustion of compact domains $K_0=
\emptyset \subset K_1 \subset K_2  \cdot \cdot \cdot  
$ for $\S_T(\a, \b)^k$, we 
then obtain a 
gluing map
$\ti \G l_k^T: U_T(\a, \b)_k \to \N_T(\a, \b)
$, where 
the domain of
gluing $U_T(p, q)_k \subset 
\S_{T}(\a, \b)_k \times 
\Bbb R^{k-1}_+$ is defined by the  condition
$$
r_i \in (0, \bar r^{-1}(K_j)]
$$
for $(u, \bold r)$ with 
$u \in K_j - K_{j-1}$. This is a smooth map.
\endproclaim
\demo{Proof} We present the case $k=2$.
The
other cases can be treated in a similar way.
Consider $u = (u^1, u^2) \in \S_T(\a, \b)_2$
and $r \in \Bbb R_+$.
For the right inverse $Q$ given
in [$\WangYe$, Lemma 4.16] we note the following
indentity
$$\|Q_{\tau_{r}(v)}\|_{\tau_r(\d_F)}
=\|Q_v\|. \tag 7.5 $$
We set $Q_1 = Q_{u^1}, Q_2= Q_{\tau_{r}(u^2)}$
and define
$\ti Q: \O^+_{1, \d}
\oplus \Ga^-_{1, \d}
\to \O^1_{2, \d}
\oplus \Ga^+_{2, \d}$
by $$\ti  Q(\zeta)=\tau_{r}(\eta_1)Q_1(\tau_{r}(
\eta_1) \zeta)+\tau_{r}(\eta_2)
Q_2(\tau_{r}(\eta_2) \zeta).$$
Then we have  for $\bar u= u \sharp R$ and $w=w_{\bar u}$
$$
\|d\hbox{\bf SW}_{\bar u}
\circ \ti Q-Id\|_{w} \leq C e^{-(\d_0-\d) r}.
$$
 Hence we obtain for large $r$ a
  right inverse $ Q = \ti Q\circ (d\hbox{\bf SW}_{\bar u}\circ
\ti Q)^{-1}$ of $d\hbox{\bf SW}_{\bar u}$.
By [$\WangYe$, Lemma 4.16],
we have  the desired norm estimate for $Q$.

By the implicit function theorem,
$\SW(\bar u+ Q) $ is a diffeomorphism from
a neighborhood $U$ of $0 \in 
\O^+_{1, \d} \oplus \Ga^-_{1, \d}$ onto
a neighborhood $\ti U$ of $\SW(\bar
u)  \in \O^+_{1, \d} \oplus \Ga^-_{1, \d}$, where 
the size of $U, \ti U$
depends only on $K$, and is independent of 
$r$.  
By Lemma 7.4,  0
must be contained in $\ti U$ if $r$
is large enough.
\qed
\enddemo

Descending to quotient, the 
  gluing maps  $\G l^T_k$ and $\ti \G l_k^T$ induce  
gluing maps  $[ \G l^T_k]$ and $[\ti \G l^T_k]$ into 
$\underline \M_T^0(S_{\a}, S_{\b})$. (For  
simplicity, we omit the indication of the time translation
quotient in the notations $[ \G l^T_k]$ and 
$[\ti \G l^T_k]$.) 
On the other hand, by a simple limit process we 
extend $\ti \G l^T_k$ to $\bar U_T(\a, \b)_k \subset 
\S_T(\a, \b)_k \times [0, \infty)^{k-1}$.
Correspondingly, we 
extend $[\ti \G l^T_k]$ to 
$\bar U_T(\a, \b)_k$. 
This extended map is a continuous map into
$\underline \MM_T^0(S_{\a}, S_{\b})$. 
It is smooth in the interior. Moreover, its restriction 
to each open  boundary stratum of $\bar U_T(\a, \b)_k$ is 
a smooth map into the corresponding stratum of 
$\underline \MM_T^0(S_{\a}, S_{\b})$.

The behavior of $\ti \G l_k^T$ along 
the boundary strata of $\bar U_T(\a, \b)_k$ 
is easy to analyse. Indeed, the relations 
given in the following lemma are easy 
consequences of the construction of 
the gluing maps.

\proclaim{Lemma 7.6} We have  
$$
\ti \G l^T_k(u, (0,...,0))=u,
\tag 7.6
$$
$$
\ti \G l^T_k((u_1,...,u_k), (0, r_2,...,r_{k-1}))
=(g(u_2,...,u_k, r_2,...,r_{k-1})^* u_1,
\tag 7.7$$
$$ 
\ti \G l^T_{k-1}((u_2,...,u_k), (r_2,...,r_{k-1}))),
$$
etc., where $g(u_2,...,u_k, r_2,...,r_{k-1}) 
\in \G(Y)$ is determined by the temporal 
transformation from $\G l_{k-1}$ to $\G l^T_{k-1}$.
\endproclaim

Before proceeding, we introduce some notations 
and terminologies. First, the space
$L_{2, \d}(\A_2(Y) \times \Ga_2(Y), 
\A_2(Y) \times \Ga_2(Y))$ (cf. [$\WangYe$, Section 6])
has a natural double affine structure given by
$$
L_{2, \d}(\A_2(Y) \times  \Ga_2(Y), \A_2(Y)
\times \Ga_2(Y)) \cong H,
$$
where 
$$H = (\oplus_{1 \leq j \leq 2}
(\O_2^1(Y) \oplus \Ga_2(Y))) \oplus 
(\O_{2, \d}^1(X) \oplus \Ga^+_{2, \d}(X))
\tag 7.8
$$
and the correspondence is given by 
$$
\eqalign {
u \to \Phi(u)= (End_-(u)-u_0, End_+(u)-u_0, u- (u_0 +
\chi (End_-(u)-u_0) +  \cr 
(1-\chi) (End_+(u)
-u_0) )\cr }
\tag 7.9
$$ 
with $\chi $ being the cut-off function introduced in 
[$\WangYe$, Section 4].  The distance $d_{2, \d}$ 
given in [$\WangYe$, Definition 4.6] can be thought of 
as induced from the natural norm on $H$ via this 
correspondence.

Next consider 
$$L_{2, \d}(\A_2(Y) \times \Ga_2(Y), 
\A_2(Y) \times \Ga_2(Y))_k \equiv \{(u_1,...,u_k) 
\in L_{2, \d}(\A_2(Y) \times \Ga_2(Y), $$
$$
\A_2(Y) \times \Ga_2(Y) )^k: End_+(u_j)=End_-(u_{j+1}) \}.
$$
 The distance $d_{2, \d}$ naturally extends to 
this space, which we still denote by $d_{2, \d}$.
On the other hand, there is a correspondence $\Phi_k$ 
similar to 
$\Phi$ with the model $H$ replaced by 
$$
H_k= (\oplus_{1 \leq j \leq k+1} (\O^1_2(Y)
\oplus \Ga_2(Y)) )  \oplus (\oplus_{1 \leq 
j\leq k  } (\O^1_{2, \d}(X) \oplus \Ga^+_{2, 
\d}(X)).
\tag 7.10
$$
The norm on $H_k$ will be denoted by $\| \cdot \|_{H_k}$.

The distance $d_{2, \d}$ descends to  $\S_T(\a, \b)_k$, which  
is a submanifold of $L_{2, \d}(\A_{2}(Y) \times \Ga_2(Y), 
\A_2(Y) \times \Ga_2(Y))_k$, and  induces its 
natural topology. We also obtain 
the tangent spaces of $\S_T(\a, \b)_k$ as suitable subspaces 
 of  $H_k$. 
 
 We have the product distance on $L_{2, \d}(\A_2(Y) \times 
 \Ga_2(Y), \A_2(Y) \times \Ga_2(Y))_k 
 \times \Bbb R^{k-1}_+$, denoted 
 again by $d_{2, \d}$, and product norm on 
 $H_k \oplus \Bbb R^{k-1}$, denoted by $\| \cdot \|_{H_k
 \oplus 
 \Bbb R^{k-1}}$.  They descend to $\S_T(\a, \b)_k 
 \times \Bbb R^{k-1}_+$ and its tangent spaces 
 respectively.

Next we introduce another norm 
which corresponds to the distance 
$d_{2, \bold r}$ introduced 
in [$\WangYe$, Definition 6.12] for measuring 
piecewise exponential convergence.

\definition{Definition 7.7}
For $\bold r \in \Bbb R_+^{k-1}$ and 
$(w_1, w_2, v)  \in T_u\S_T(\a, \b)$ we set
$$ 
\|(w_1, w_2, v)\|_{\bold r} ^2= \sum_{1\leq j \leq 2}
\|w_i\|_2^2 + \sum_{1\leq j \leq k-1} \| v(\cdot, 
r_i)\|_2^2 
+ \|v-Int^*_{\bold r} (v)\|_{2, w_{\bold r}}^2,
\tag 7.11
$$
where
$$
Int^*_{\bold r}(v)=  
\tau_{r_1}(\eta) v(\cdot, r_1) 
+...+ $$
$$\tau_{2r_1+...2r_{k-2}+r_{k-1}}(\eta) 
v(\cdot, 2r_1+...2r_{k-2}+ r_{k-1})
$$
and $\eta$ is the cut-off function introduced in 
[$\WangYe$, Section 6].
\enddefinition

\proclaim{Lemma 7.8} $[\ti \G l^T_k] $  is a local 
diffeomorphism along $U_T(\a, \b)_k$. 
Equivalently, $[ \G l^T_k] $ is a local 
diffeomorphism along the domain corresponding 
to $U_T(\a, \b)_k$ (the correspondence is 
in terms of the map $Inv$).
Indeed,
for each compact set $K$ in $\S_T(\a, \b)_k$,
there are  positive numbers  $r_K, R_K$ and 
$\rho_K$ with the 
following property. Let $u \in K$ and $ 
\bold r= (r_1, ..., r_{k-1}) 
\in \Bbb R_+^{k-1}$ with $ r_i \in (r_K, \infty),
i=1,...,k-1$. Then there is a neighborhood 
$U(u, \bold r)$ of $(u, \bold r)$
such that the restriction of $[ \G l^T_k]$ 
to it is a diffeomorphism onto the 
distance ball $B_{\rho_K}([\G l^T_k](u, 
\bold r))$, where the distance 
is the $(2, \bold r)$-distance $d_{2, \bold r}$. Moreover, $U(u, \bold r)$ contains
$B_{r_K}(u) \times B_{R_K}(\bold r)$, where 
the distance is $d_{2, \d}$. 
\endproclaim
\demo{Proof} Since $\S_T(\a, \b)$ is 
a presentation model for $\underline \M_T^0(S_{\a},
S_{\b})$, we obtain from $[\G l^T_k]$ an induced 
gluing map $ \G l^0_k: U_T(\a, 
\b)_k \to \S_T(\a, \b)$. We show that 
$\G l^0_k $ is a local diffeomorphism with 
size control corresponding to the size control 
stated in the lemma.

Fix a compact domain $K$ in $\S_T(\a, \b)_k$, e.g. $K=K_j$
 for some 
$j$.  Consider $u \in K$ and $\bold 
r \in \Bbb R_+^{k-1}$ with $r_i \in
[2\bar r_K, \infty)$. We write 
$u^*= \G l^0_k(u, \bold r)$. 
Solving the equation $\SW(\Phi_k^{-1}(u + v+ v^{\perp}))=0$
(this means $\SW=0 $ for each portion in $u+v+ v^{\perp}$)
for $v \in T_u\S_T(\a, \b)_k,
v^{\perp} \in T_u\S_T(\a, \b)_k$ by using the implicit function 
theorem and elliptic estimates, we obtain 
a coodinate map $\Phi_u: B_{r_K}(0) \to 
U(u)$ with smooth control,
such that $U(u)$ contains the ball
$B_{\rho_K}(u)$. Similarly, we obtain 
a coordinate map $\Psi_{u^*}: B^{\bold r}_{r_K^*}(0)
\to U(u^*)$ with smooth control,
where $U(u^*)$
contains the ball $B_{r^*_K}^{\bold r}(u^*)$.
Here, the superscript $\bold r$ means that 
the ball is a distance ball with respect to 
the norm $\| \cdot \|_{\bold r}$ or 
the distance $d_{2, \bold r }$.

By the construction of $\G l^0_k$ and 
elliptic estimates, we can 
choose the numbers $ r_K, r_K^*, \rho_K$
and $
\rho_K^*$ suitably, such that 
the composition 
map 
$$\G l^u \equiv \Psi_{u^*}^{-1} \cdot \G l^0_k 
 \cdot  (\Phi_{u}
\times Id)$$ 
goes from $B_{r_K} \times 
 B_{\bar r_K /2}(\bold r)$ 
into $B_{r_K^*}^{\bold r}(0)$.

Now we proceed to compute the differential of 
$\G l^u$. For this purpose, we introduce  
another gluing map.  For $v \in B_{r_K}(0)$
and $\bold r' \in B_{r_K}(\bold r)$,
we consider the 
equation
$$
\SW(\#(\Phi_k^{-1}(v), \bold r') 
+ Q v')=0.
$$
By the proof of Lemma 7.4 we obtain solutions 
$$
v' = \zeta^*(v, \bold r')
$$
for $v, \bold r'$ in balls of size depending 
only on $K$.
We set 
$$
\sigma(v, \bold r') = \#(\Phi_k^{-1}(v), \bold r')
+ Q \zeta^*(v, \bold r').
\tag 7.12
$$
Converting $\sigma$ into the slice $\S_T(\a, \b)$
by a transformation $P$ involving the time translation and 
gauge transformations, we obtain a new map which 
we denote by $\G l^{u}_0$.
By tangency,
the differential $d \G l^u_0|_u$ equals 
the differential $d \G l^u|_u$. Since the transformation
$P$ can easily be controlled, it suffices for 
the purpose of computing $d \G l^u|_u$ to compute 
 $d \sigma|_u$. 

Consider e.g. the case $k=2$. The 
general case is similar. 
Consider the above $u=(u_1, u_2)$ and $\bold r =r$.
The contribution of $Q \zeta^*$ to $\Sigma$ is small 
and can be absorbed into that of  
the first term.  Thus, we have 
$$
d \sigma|_{(u, \bold r)}(0, (v_1, v_2), 0) \cong
(v_1, 0) \sharp r  + (0, v_2) \sharp r.
\tag 7.13
$$

Using the exponential decay properties 
of the tangent vectors $v=(v_1, v_2)$
and the finite dimensionality of the tangent 
space
we deduce 
$$\|d \sigma|_{(u, \bold r)} (0, (v_1, v_2), 0)\|_{\bold r} 
\geq  c(K) \|(v_1, v_2)\|_{2, \d}
\tag 7.14
$$
for a positive constant $c(K)$.
Here the finite dimensionality is used in the 
following way: we consider an orthonormal base 
of tangent vectors and derive the estimate (7.14) 
for them first. Then (7.14) follows for general 
tangent vectors $(0, v_1, v_2, 0)$.

On the other hand, we have
$$
d \sigma|_{(u, \bold r)} (0, 0, 1)
\cong \tau_r (\frac{\partial u_2}{\partial t}).
\tag 7.15
$$
Using the arguments in the proof of Lemma 7.2 and compactness 
we derive a positive lower bound for 
$\int_{t=0} |\pa u_2 / \pa t|^2$.  Consequently,
we deduce   a positive 
lower bound on  
$\|d \Sigma|_{(u, \bold r)} (0, 0, 1)\|_{2, \bold r}.$ 
Combining this with  (7.11) and applying the 
transversality of the model $\S_T(\a, \b)$ 
with respect to the time translation action,
we then deduce 
$$
\| d \Sigma|_{(u, \bold r)} (0, (v_1, v_2), r')\|_{\bold r}
\geq c(K) (\|(v_1, v_2)\|_{2, \d}^2+ (r')^2)^{1/2}.
\tag 7.16
$$
Finally, we have
$$
d \sigma|_{(u, \bold r)} (w_1, w_2, w_3, 0, 0, 0) \cong 
\chi w_1 + \tau_r(\eta) w_2 + (1-\chi) w_3.
\tag 7.17
$$
Combining the  above estimates and arguments
we arrive at
$$
\|d \Sigma|_{(u, \bold r)}\| \geq c(K) >0,
\tag 7.18
$$
where the operator norm is defined with 
respect to the norms $\| \cdot \|_{H_k \oplus \Bbb R^{k-1}_+}$
and $\|\cdot \|_{\bold r}$.
Note that it is easy to handle the cross interactions 
between   $(w_1, w_2, w_3)$ and 
$(v_1, v_2, r')$ here.
The estimate (7.18) immediately leads to  
$$
\|d  \G l^u|_{(u, \bold r)}\| \geq c(K) >0.
\tag 7.19
$$
We can compute the differential $d \G l^u$ at other 
points in  a similar way
and derive an upper bound on its norm.  
It involves only minor modifications of the 
above computations. 

By the implicit function theorem, $\G l^u  $ 
is a diffeomorphism from a neighborhood of 
$(0, \bold r, 0)$ onto a neighborhood of 
$0$ with 
the desired size control. This implies that 
$\G l^0_k$ is a diffeomorphism from 
a neighborhood of $(u, \bold r)$ onto a 
neighborhood of $u^*$ with the desired 
size control. The desired properties of $[\G l^T_k]$
follow. 
\qed 
\enddemo

\proclaim{Lemma 7.9} The restriction of 
$[\ti \G l_k^T]$ to each open stratum of 
$\bar U_T(\a, \b)_k$ is a local diffeomorphism  
into the corresponding open stratum of 
$\underline \MM_T^0(S_{\a}, S_{\b})$ with 
size control similar to that in Lemma 7.8.
\endproclaim

The poof of this lemma is similar to that 
of Lemma 7.8, we omit the details.

\proclaim{Lemma 7.10} 
We can choose the domain 
$U_T(\a, \b)_k$ suitably such that 
$[\ti \G l^T_k]$ is a homeomorphism from 
$\bar U_T(\a, \b)_k$ onto a 
neighborhood of $\underline \M_T^0(S_{\a},
S_{\b})_k$ in the moduli 
space $\underline \MM^0_T(S_{\a},
S_{\b})$.
\endproclaim
\demo{Proof} First we show that the image of 
$\bar U_T(\a, \b)_k$ is a neigborhood of 
$\underline \M_T^0(S_{\a}, S_{\b})_k$ in 
$\underline \MM_T^0(S_{\a}, S_{\b})$.  Assume the contrary.
Then there would be a sequence $\omega_j$ 
in $\underline \MM_T^0(S_{\a}, S_{\b})$ converging 
to some $ \omega \in \underline \M_T^0(S_{\a},
S_{\b})_k$ such that none of $\omega_j$ is 
in the image of $[\ti  \G l^T_k]$.  We consider the 
case $\omega_j \in \underline \M^0_T(S_{\a},
S_{\b})$, while the other cases are 
similar.  
We have representatives $u_j \in 
\omega_j, u \in \omega$ such that 
$u_j \in \S_T(\a, \b),  u \in 
\S_T(\a, \b)_k$.  By Lemma 7.3, $\tau_{-t_k}(u_j)$ 
converge to 
$u$  piecewise exponentially, where $t_k$ is 
a suitable sequence of numbers.  Hence 
there is a sequence $\bold r_j =
(r_{k, 1},...,r_{j, k-1})$ such that   
$$d_{2, \bold r_j }(\tau_{-t_k}(u_j), u \sharp \bold r_j)
\to 0.
$$
It follows then that 
$$
d_{2, \bold r_j}(\omega_j,  [\G l^T_k](u, \bold r_j))
\to 0.
$$
Applying Lemma 7.8 we deduce that for large 
$j$, $\omega_j$ 
is in the image of $[\G l^T_k]$, and hence in
the image of $[\ti \G l]_k$. This is a 
contradiction.

Next we show that $[\ti \G l^T_k]$ is injective.
First note that the images of different strata 
of $\bar U_T(\a, \b)_k$
under $[\ti \G l^T_k]$ are disjoint.  Hence we can 
consider each stratum individually.  Consider 
e.g. the top stratum, which is contained in 
$\S_T(\a, \b)_k \times \Bbb R^{k-1}_+$. Lower 
strata can be handled in a similar way, where 
we apply Lemma 7.9 
instead of Lemma 7.8. 
Assume   $[\ti \G l^T_k] (u, \bold r) = 
[\ti \G l^T_k](\ti u,  \ti {\bold r})$ for $u, \ti u\in
\S_T(\a, \b)_k, \bold r, \ti {\bold r} \in
\Bbb R^{k-1}_+.$  Consider the 
path $\theta(t)= (u(t), 
\bold r(t)) = (u, t \bold r), t \in 
[0, 1]$. By Lemma 7.8 we can find a 
path $\ti \theta(t)= 
(\ti u(t), \ti {\bold r}(t)),  t \in (0, 1]$ such that 
$ \ti \theta(1) = (\ti u, \ti {\bold r})$ and 
$[\ti \G l^T_k] (\ti \theta(t))= [\ti \G l^T_k](\theta(t)).$
By  the identity (7.6) in
Lemma 7.6, we deduce that 
the path $\ti \theta$ extends continuously 
to $\ti \theta(0)=(u, (0,...,0))$. 

We claim that for $t$ small enough, $\ti \theta(t)=
\theta(t)$. For simplicity, we consider the 
case $k=2$. The general case is similar. 
In this case we have  $u=(u_1, u_2), \bold r =r, 
\ti u= (\ti u_1, \ti u_2), \ti {\bold r} = \ti r$ etc.

There is a unique number $R^+_{u_1}$ such 
that $E(u_1, Y \times [R^+_{u_1}, \infty)) 
= E_0/8$ and $E(u_1, Y \times (-\infty, R^+_{u_1}])
\geq 7E_0/8.$ Similarly, there is a unique 
number $R^-_{u_2}$ such that 
$E(u_2, Y\times (-\infty, R^-_{u_2}])
=E_0/8$ and $ E(u_2, [R^-_{u_2}, \infty)) 
\geq 7E_0/8.$ It follows that $u
\# (tr)^{-1}$ has an ``energy valley" of length 
$2(tr)^{-1} -R^+_{u_1}-R^-_{u_2} + o(1)$ with 
$o(1) \to 0$ as $r \to 0$. Here an
``energy valley" of $u \#(tr)^{-1}$  means an interval $I=
[R_1, R_2]$ such that 
$$
E(u \# (tr)^{-1}, Y \times I) = E_0/4,
$$
$$
E(u \# (tr)^{-1}, Y \times (-\infty, 
R_1] ) \geq  3E_0/4, E(u 
\# r, Y \times [R_2, \infty)) \geq 3E_0/4.
$$
By the gluing construction and the 
consequent elliptic estimates, $ \G l_k(u, (tr)^{-1})$
also has an energy valley of length 
$
2(tr)^{-1} -R^+_{u_1}-R^-_{u_2} +o(1) $.

Similarly, $ \G l_k(\ti u(t), \ti r(t)^{-1})$ has 
an energy valley of length $2 \ti r(t)^{-1}
-R^+_{\ti u_1(t)} -
R^-_{\ti u_2(t)} +o(1)$. 

We deduce
$$
2 \ti r(t)^{-1}
-R^+_{\ti u_1(t)} -
R^-_{\ti u_2(t)} +o(1)=
2(tr)^{-1} -R^+_{u_1}-R^-_{u_2} +o(1).
\tag 7.15
$$

 Since $\ti u(t) 
\to u$ as $t \to 0$, we have 
$$
R^+_{\ti u_1(t)} \to R^+_{u_1}, R^-_{\ti 
u_2(t) } \to R^-_{u_2}.
\tag 7.16
$$
Combining (7.15) with (7.16) we infer that 
$$\ti r(t)^{-1}= (tr)^{-1}+ o(1).$$

By Lemma 7.8 we then deduce that 
$u(t)= \ti u(t), \ti r(t)= tr$ for small $t$.  Now we 
can repeat the argument with a small, 
positive $t_0$ replacing 
$t=0$. By extension and continuity, we 
conclude that $\ti u(1) = u(1), \ti r(1) = r$, i.e.
$\ti u =u, \ti r = r$.

Finally, Propostion 6.11  implies that 
$[\ti \G l_k^T]$ is proper. Consequently, it 
is a homeomorphism.
\qed
\enddemo

\demo{Proof of Proposition 7.1}
We define $\bold F_k$ to be $[\ti \G l^T_k]$. 
The previous lemmas immediately lead to the desired 
properties of $\bold F_k$. Note in particular 
that the smooth transition property (property (4))
is an easy  consequence of Lemma 7.7 and Lemma 7.8.
\qed
\enddemo

\definition{Remark 7.12} All the above arguments 
extend 
straightforwardly to moduli spaces of 
transition trajectories as can easily be seen.
Note that constant trajectories appear in the 
invariance proof. They have zero energy, and hence 
the energy valley argument in the proof of Lemma 7.10
does not apply to them. However, their moduli 
spaces are smoothly compact, namely the 
piecewise exponential convergence here is 
actually exponential convergence to smooth 
trajectories. Hence  there is 
no need of using the energy valley argument for them.
\enddefinition

Finally, we sketch the gluing process in 
terms of the fixed-end model.    First,
we need a presentation model for 
$\underline \MM^0(p, q; \S \W_0)$, where 
the underline means quotient under the 
twisted time translation action introduced 
in [$\WangYe$, Definition 6.16]. The twisted 
time translation action can be handled 
in a similar way to the time translation action
for the temporal model. On the other hand,
we can use the global Columb gauge to 
handle the based gauge action. Namely we 
consider the slice given by 
$d_{\d}^*(A-A_0)=0$. Here there is 
a delicate point when $p$ or $q$ is 
reducible. In this case, additional 
gauges in the isotropy 
group $G_p$ or $G_q$ (see 
[$\WangYe$, Section 4]) are involved.  As a consequence,
we can only obtain 
local slices instead of global slices. 
The remaining parts of the process are 
similar to the case of the temporal model.
We omit the details. 
In comparison, we see two advantages of 
the temporal model approach. The first 
is the fact that it is canonical, while 
the fixed-end model involves the choice of 
$\S \W_0$. The second is that we have 
global slices and presentation models 
for the temporal model, while we only 
have local slices and presentation models 
for the fixed-end model.


\Refs

\refstyle{1}
\widestnumber\key{JPW}
 
\ref\key 1
\by M.F. Atiyah
\paper New invariants of 3- and 4-dimensional manifolds
\jour Symp. Pure Math.
\vol 48
\yr 1988
\pages 285-299
\endref

\ref\key 2
\by M.F. Atiyah and R. Bott
\paper The moment map and equivariant cohomology
\jour Topology
\vol 23
\yr 1984
\pages 1-28
\endref

\ref\key 3
\by M.F.Atiyah, V.K. Patodi and I.M. Singer
\paper ~Spectral asymmetry and Riemannian geometry I, II, III
\jour Math. Proc. Cambridge. Soc.
\vol 77 
\yr 
\pages 
\endref


\ref\key  4
\by D. M. Austin and P. J. Braam
\paper Equivariant Floer theory and gluing 
Donaldson polynomials 
\jour Topology
\vol 35 
\yr 1996
\pages 167-200
\endref

\ref\key 5
\by D. M. Austin and P. J. Braam
\paper Morse-Bott theory and equivariant 
cohomology
\inbook Floer Memorial Volume
\bookinfo  Birkh{\"a}user 
\yr 1996
\pages 123-164
\endref

\ref\key 6
\by N. Berline, E. Getzler and M. Vergne
\book Heat kernels and Dirac operators
\bookinfo Springer-Verlag
\yr 1992
\endref

\ref\key 7
\by    A. Borel
\paper Seminar on Transformation groups
\inbook Annals of Mathematical Studies 46
\bookinfo  Princeton University Press 
\yr 1960
\pages 
\endref

\ref\key 8
\by   H. Cartan
\paper  La transgression dans un groupe de Lie
et dans un espace fibr{\'e} principal
\inbook Colloque de Topologie (Espaces Fibr{\'e})
\bookinfo  C. B. R. M. Bruxells 
\yr 
\pages 57-71
\endref


\ref\key 9
\by W. Chen
\paper Casson's invariant and Seiberg-Witten 
gauge theory 
\inbook preprint
\bookinfo 
\yr 
\pages 
\endref

\ref\key 10
\by R. Cohen, J. Jones and G. Segal
\paper Floer's infinite dimensional Morse theory and homotopy theory
\inbook Floer Memorial Volume
\bookinfo Birkh{\"a}user
\yr 1996
\pages 297-326
\endref

 

\ref\key 11
\by S.K. Donaldson
\paper The orientation of Yang-Mills moduli spaces and 4-manifold topology
\jour J. Diff. Geom.
\vol 26
\yr 1987
\pages 397-428
\endref




\ref\key 12 
\bysame 
\paper  The Seiberg-Witten 
equations and the 4-manifold 
topology
\jour Bull. Amer. Math. Soc.
\vol
\yr 1996
\pages 
\endref


\ref\key 13
\by A. Floer
\paper Holomorphic spheres and symplectic 
fixed points
\jour Comm. Math. Phys.
\vol
\yr
\pages 
\endref

\ref\key 14 
\by A. Floer
\paper An instanton invariant for 3-manifolds
\jour Comm. Math. Pyhs.
\vol 118
\yr 1988
\pages 215-240
\endref




\ref\key 15 
\by Fukaya
\paper Floer homology of connected sum of homology 3-sphere
\jour Topology
\vol 36
\yr  1996
\pages 
\endref






\ref\key 16
\by P. B. Kronheimer and T. S. Mrowka
\paper The genus of embedded surfaces in the projective space
\jour Math. Res. Letters
\vol 1
\yr 1994
\pages 797--808
\endref



\ref\key 17
\by  H. B. Lawson and M.-L Michelsohn
\book  Spin geometry
\bookinfo ~Princeton, New Jersey
\yr 1989
\endref

\ref\key 18
\by Y. Lim
\paper The equivalence of Seiberg-Witten and 
Casson invariants for homology 3-spheres
\jour preprint
\yr
\endref

\ref\key 19
\by M. Marcolli and L. Wang
\paper Equivariant Seiberg-Witten Floer homology
\jour preprint
\yr
\endref






\ref\key 20
\by D. Salamon and E. Zehnder
\paper Morse theory for periodic solutions of Hamiltonian system and
the Maslov index
\jour Comm. Pure Appl. Math.
\vol 45
\yr 1992
\pages 1303-1360
\endref


\ref\key 21 
\by C. Taubes 
\paper   Casson's invariant and gauge theory 
\jour J. Diff. Geom.
\vol 31
\yr 1990
\pages 547-599
\endref

\ref\key 22 
\bysame
\paper The Seiberg-Witten invariants and the Gromov invariants
\jour preprint 
\yr
\pages
\endref





\ref\key 23
\by Guofang Wang and Rugang Ye 
\paper Bott-type and equivariant Seiberg-Witten
Floer homology I 
\jour DG9701010 
\vol  
\yr 1997
\pages 
\endref

\ref\key 24
\by Guofang Wang and Rugang Ye 
\paper Equivariant and Bott-type Seiberg-Witten
Floer homology: Part I 
\jour math.GT/9901058
\vol 
\yr 
\pages 
\endref

\ref\key 25
\by E. Witten
\paper Monopoles and 4-manifolds
\jour Math. Res. Letters
\vol 1
\yr 1994
\pages 769--796
\endref


\ref\key 26
\by Rugang Ye
\paper Equivariant and Bott-type Seiberg-Witten
Floer homology: Part III
\jour in preparation
\vol
\yr
\pages
\endref


\ref\key 27
\by Rugang Ye
\paper Gromov compactness theorem for pseudo
holomorphic curves
\jour Tran. Amer. Math. Soc. 
\vol
\yr 1994
\pages
\endref

\endRefs
\enddocument